\documentclass[11pt]{article}

\parskip=\medskipamount
\def\fpd#1#2{{\displaystyle\frac{\partial #1}{\partial #2}}}
\def\spd#1#2#3{{\displaystyle\frac{\partial^2 #1}
{\partial #2\partial #3}}}
\def\lie#1{{\cal L}_{#1}}
\def\vf#1{\frac{\partial}{\partial #1}}
\def\ad{\mathop{\mathrm{ad}}\nolimits}

\def\conn#1#2#3{\setbox1=\hbox{$\scriptstyle{#2}{#3}$}%
\setbox2=\hbox to\wd1{$\hfil\scriptstyle{#1}\hfil$}
\Gamma^{\!\box2}_{\!\box1}}

\def\clift#1{#1^{\scriptscriptstyle{\mathrm{C}}}}
\def\hlift#1{#1^{\scriptscriptstyle{\mathrm{H}}}}
\def\vlift#1{#1^{\scriptscriptstyle{\mathrm{V}}}}
\def\mech#1{#1^{\scriptstyle{\mathrm{m}}}}
\def\submech#1{#1_{\scriptstyle{\mathrm{m}}}}

\def\A{\mathcal{A}}

\def\R{\mathcal{R}}

\def\onehalf{{\textstyle\frac12}}

\font\frak=eufm10 scaled\magstep1
\def\goth #1{\hbox{{\frak #1}}}
\def\g{\goth{g}}

\def\J{N_\mu}
\def\Hb{\hlift{{\bar E}}}
\def\Hh{\hlift{{\hat E}}}
\def\F{\hlift{{\check E}}}

\newcommand{\mybox}[1]{\makebox(0,0){\footnotesize{#1}}}
\def\vectorfields#1{{\cal X}(#1)}

\def\pM{{\pi}^{\scriptscriptstyle M}}
\def\psiM{\psi^{\scriptscriptstyle M}}
\def\psiN{\psi^{\scriptscriptstyle N}}
\def\psiTM{\psi^{\scriptscriptstyle TM}}
\def\Ad{\mathop{\mathrm{ad}}\nolimits}
\def\ad{\Ad}

\def\x{y}
\def\y{z}
\def\z{x}

\begin{document}

\title{Routh's procedure for non-Abelian symmetry groups}

\author{M.\ Crampin${}^{a}$ and T.\ Mestdag${}^{a,b}$\\[2mm]
{\small ${}^a$Department of Mathematical Physics and Astronomy, Ghent University}\\
{\small Krijgslaan 281, B-9000 Ghent, Belgium}\\[1mm]
{\small ${}^b$ Department of Mathematics, University of Michigan}\\
{\small 530 Church Street, Ann Arbor, MI 48109, USA}}

\date{}

\maketitle

{\small {\bf Abstract.} We extend Routh's reduction procedure to an
arbitrary Lagrangian system (that is, one whose Lagrangian is not
necessarily the difference of kinetic and potential energies) with a
symmetry group which is not necessarily Abelian.  To do so we
analyse the restriction of the Euler-Lagrange field to a level set
of momentum in velocity phase space.  We present a new method of
analysis based on the use of quasi-velocities.  We discuss the
reconstruction of solutions of the full Euler-Lagrange equations
from those of the reduced
equations.\\[2mm]
}

\section{Introduction}

Routh's procedure, in its original form (as described in his
treatise \cite{Routh}), was concerned with eliminating from a
Lagrangian problem the generalized velocities corresponding to
so-called ignorable or cyclic coordinates.  Let $L$ be a Lagrangian
on ${\bf R}^n$ that does not explicitly depend on $m$ of its base
variables, say the coordinates $\theta^a$. From the Euler-Lagrange
equations for these coordinates,
\[
\frac{d}{dt}\left(\fpd{L}{{\dot\theta}^a}\right)-\fpd{L}{\theta^a}=0,
\]
we can immediately conclude that the functions $\partial
L/\partial\dot{\theta}^a$ are constants, say
\[
\fpd{L}{\dot{\theta}^a}=\pi_a;
\]
these equations express the conservation of generalized momentum.
Routh's idea is to solve these equations for the variables
${\dot\theta}^a$ and to introduce what he calls the `modified
Lagrangian function', the restriction of the function
\[
L'=L-\fpd{L}{{\dot\theta}^a}{\dot\theta}^a
\]
to the level set where the momentum is $\pi_a$.  One can easily verify
that the $(n-m)$ Euler-Lagrange equations for the remaining variables
$x^i$ can be rewritten as
\[
\frac{d}{dt}\left(\fpd{L}{{\dot x}^i}\right)-\fpd{L}{x^i}=0 \quad
\Rightarrow \quad \frac{d}{dt}\left(\fpd{L'}{{\dot
x}^i}\right)-\fpd{L'}{x^i}= 0.
\]
For example, if the Lagrangian takes the form
\[
L(x,\theta,\dot x,\dot\theta)= \onehalf k_{ij}(x){\dot x}^i{\dot
x}^j + k_{ia}(x) {\dot x}^i{\dot \theta}^a + \onehalf
k_{ab}(x){\dot \theta}^a{\dot \theta}^b-V(x),
\]
the conservation of momentum equations read $k_{ia}{\dot
x}^i+k_{ab}{\dot\theta}^b=\pi_a$, and they can be solved for the
variables ${\dot\theta}^a$ if $(k_{ab})$ is a non-singular matrix.
The modified Lagrangian function is
\[
L'(x,{\dot x})=
\onehalf(k_{ij}-k^{ab}k_{ia}k_{jb}){\dot x}^i{\dot x}^j
+k^{ab}k_{ia}\pi_b{\dot x}^i -
(V+\onehalf k^{ab}\pi_a\pi_b),
\]
where $k^{ab}$ denotes a component of the matrix inverse to
$(k_{ab})$, in the usual way. Clearly, the advantage of this
technique is that the reduced equations in $L'$ involve only the
unknowns $x^i$ and $\dot{x}^i$; they can in principle be directly
solved for the $x^i$, and the $\theta^a$ may then be found (if
required) from the momentum equation.

A modern geometric interpretation of this reduction procedure can be
found in e.g.\ \cite{MRbook}.  The above Lagrangian $L$ is of the form
$T-V$, where the kinetic energy part is derived from a Riemannian
metric (i.e.\ we are dealing with a so-called simple mechanical
system).  The function $L$ is defined on the tangent manifold of a
manifold of the form $M = S\times G$ (in this case ${\bf R}^n$) and it
is invariant under an Abelian Lie group $G$ (in this case the group of
translations ${\bf R}^m$).  The main feature of the procedure is that the
modified Lagrangian function and its equations can be defined in terms
of the coordinates on $S$ only.  However, to give the definition of
the modified function an intrinsic meaning, we should define this
function, from now on called the Routhian, rather as the restriction
to a level set of momentum of
\[
\R=L-\fpd{L}{{\dot\theta}^a}({\dot\theta}^a
+\Lambda^a_i {\dot x}^i),
\]
with $\Lambda^a_i= k^{ab}k_{ib}$, i.e.\
\[
\R(x,{\dot x})= \onehalf(k_{ij}-
k^{ab}k_{ia}k_{jb}){\dot x}^i{\dot x}^j - (V+\onehalf k^{ab}\pi_a\pi_b).
\]
The coefficients $\Lambda^a_i$ form a connection on the trivial
principal bundle $M = S\times G \to S$, usually called the
mechanical connection, and ${\dot\theta}^a +\Lambda^a_i {\dot x}^i$
is in fact the vertical projection of the vector $({\dot
x}^i,{\dot\theta}^a)$. The $(n-m)$ Euler-Lagrange equations in $x^i$
then become
\[
\frac{d}{dt}\left(\fpd{\R}{{\dot x}^i}\right)-\fpd{\R}{x^i}= -
B^a_{ij}\pi_a{\dot x}^j,
\]
where in the term on the right-hand side
\[
B^a_{ij}=\fpd{\Lambda^a_i}{x^j}- \fpd{\Lambda^a_j}{x^i}
\]
has a coordinate-free interpretation as the curvature of the
connection.

In \cite{MS,MRS}, Marsden et al.\ extended the above procedure to
the case of simple mechanical systems with a non-Abelian symmetry
group $G$ and where the base manifold has a principal bundle
structure $M\to M/G$.  The procedure has recently been further
extended to cover Lagrangian systems in general by Castrillon-Lopez
\cite{Marco}.

The most important contribution of our paper lies in the geometric
formalism we will adopt.  The bulk of the literature dealing with
different types of reduction of Lagrangian systems has relied
heavily on methods coming from the calculus of variations.  In fact,
as in e.g.\ \cite{Marco,JM,MRS}, the reduced equations of
motion are usually obtained by considering some reduced version of
Hamilton's principle.  Our method is different from those of other
authors in that it doesn't involve consideration of variations.  It
is distinctively Lagrangian (as opposed to Hamiltonian), and is
based on the geometrical analysis of regular Lagrangian systems,
where solutions of the Euler-Lagrange equations are interpreted as
integral curves of an associated second-order differential equation
field on the velocity phase space, that is, the tangent manifold of
the configuration space.  Consequently our derivation of Routh's
equations is relatively straightforward and is a natural extension
of that used by Routh in the classical case.  In particular, we will
show how Routh's equations can be derived directly from the
Euler-Lagrange equations by choosing a suitable adapted frame, or
equivalently by employing well-chosen quasi-velocities.  This line
of thinking has already provided some new insights into e.g.\ the
geometry of second-order differential systems with symmetry
\cite{Paper1}.

We deal from the beginning with arbitrary Lagrangians, i.e.\
Lagrangians not necessarily of the form $T-V$.

As in \cite{MRS}, we explain how solutions of the Euler-Lagrange
equations with a fixed momentum can be reconstructed from solutions
of the reduced equations. The method relies on the availability of a
principal connection on an appropriate principal fibre bundle. We
will introduce in fact two connections that serve the same purpose.

We describe the basic features of our approach in Section~2.  The
reduction of a Lagrangian system to a level set of momentum is
discussed in Section~3, and our generalization of Routh's procedure
is explained there.  Section~4 contains some general remarks about
using a principal connection to reconstruct an integral curve of a
dynamical vector field from one of a reduction of it.  In Section~5
we describe the two principal connections that can be used in the
specific reconstruction problem we are concerned with, while in
Section 6 we carry out the reduction in detail, first in the Abelian
case, then in general.  In Section~7 we specialize to simple
mechanical systems, in order to compare our results with those
published elsewhere.  We conclude the paper with a couple of
illustrative examples.

\section{Preliminaries}

We will be concerned with Lagrangian systems admitting non-Abelian
(that is to say, not necessarily Abelian) symmetry groups. We begin
by explaining what assumptions we make about the action of a
symmetry group.

We will suppose that $\psiM: G\times M \to M$ is a free and proper
left action of a connected Lie group $G$ on a manifold $M$.  It
should be noticed from the outset that this convention differs from
the one in e.g.\ \cite{Paper1,KN}, but resembles the one taken in
e.g.\ \cite{MRbook,MRS}.

With such an action, $M$ is a principal fibre bundle with group $G$;
we write $M/G$ for the base manifold and $\pM:M\to M/G$ for the
projection.  We denote by $\g$ the Lie algebra of $G$.  For any
$\xi\in\g$, $\tilde\xi$ will denote the corresponding fundamental
vector field on $M$, that is, the infinitesimal generator of the
1-parameter group $\psiM_{\exp(t\xi)}$ of transformations of $M$.
The Lie bracket of two fundamental
vector fields satisfies
$[\tilde\xi,\tilde\eta]=-\widetilde{[\xi,\eta]}$ (see e.g.\
\cite{MRbook}). Since $G$ is connected, a tensor field on $M$ is
invariant under the action of $G$ if and only if its Lie derivatives
by all fundamental vector fields vanish. In particular, a vector
field $X$ on $M$ is invariant if and only if $[\tilde\xi,X]=0$ for
all $\xi\in\g$.  We will usually work with a fixed basis for $\g$,
which we denote by $\{E_a\}$; then for $X$ to be invariant it is
enough that $[\tilde{E}_a,X]=0$, $a=1,2,\ldots,\dim(\g)$.

We suppose that we have at our disposal a principal connection on
$M$. For the most part it will be convenient to work with
connections in the following way.  A connection is a left splitting
of the short exact sequence
\[
0 \to M\times \g \to TM \to (\pM)^*T(M/G) \to 0
\]
of vector bundles over $M$; we identify $M\times\g$ with the
vertical sub-bundle of $TM\to M$ by $(m,\xi)\mapsto\tilde{\xi}|_m$.
Thus we may think of a connection as a type $(1,1)$ tensor field
$\omega$ on $M$ which is a projection map on each tangent space,
with image the tangent to the fibre of $\pM$.  The connection is
principal just when $\omega$ is invariant, that is, when
$\lie{\tilde{\xi}}\omega=0$ for all $\xi\in\g$.  The kernel
distribution of $\omega$ is the horizontal distribution of the
connection.  An alternative test for invariance of the connection is
that its horizontal distribution should be invariant (as a
distribution); that is, for any horizontal vector field $X$,
$[\tilde{\xi},X]$ is also horizontal for all $\xi$.  We will often
refer to a connection by the symbol of the corresponding tensor
field.

Let $\{X_i\}$ be a set of local vector fields on $M$ which are
linearly independent, horizontal with respect to $\omega$ and
invariant.  Such a set of vector fields consists of the horizontal
lifts of a local basis of vector fields on $M/G$, and in particular
we may take for the $X_i$ the horizontal lifts of coordinate fields
on $M/G$.  We then have a local basis $\{X_i,\tilde{E}_a\}$ of
vector fields on $M$.  We will very often work with such a basis,
which we call a standard basis.  The Lie brackets of pairs of vector
fields in a standard basis are
\[
[X_i,X_j] = R^a_{ij} {\tilde E}_a,\quad [X_i,{\tilde E}_a]=0,\quad
[{\tilde E}_a,{\tilde E}_b]=-C_{ab}^c {\tilde E}_c.
\]
The $R^a_{ij}$ are the components of the curvature of $\omega$,
regarded as a $\g$-valued tensor field.  The second relation simply
expresses the invariance of the $X_i$.  In the third expression the
$C_{ab}^c$ are structure constants of $\g$ with respect to the
chosen basis.

It will sometimes be convenient to have also a basis $\{X_i,{\hat
E}_a\}$ that consists entirely of invariant vector fields.  Let
$U\subset M/G$ be an open set over which $M$ is locally trivial.  The
projection $\pM$ is locally given by projection onto the first factor
in $U\times G \to U$, and the (left) action by $\psiM_g(x,h)=(x,gh)$.
The vector fields on $M$ defined by
\[
{\hat E}_a: (x,g) \mapsto \widetilde{(\Ad_{g} E_a)}(x,g) = \psiTM_g
\big({\tilde E}_a  (x,e)\big).
\]
(where $e$ is the identity of $G$) are invariant.  The relation
between the sets $\{{\hat E}_a\}$ and $\{{\tilde E}_a\}$ can be
expressed as ${\hat E}_a(x,g)= {\mathcal A}_a^b(g){\tilde E}_b(x,g)$
where $({\mathcal A}_a^b(g))$ is the matrix representing $\Ad_{g}$
with respect to the basis $\{E_a\}$ of $\g$.  In particular,
${\mathcal A}^b_a(e)=\delta^b_a$.  Since $[{\tilde E}_a, {\hat E}_b]
= 0$, the coefficients ${\mathcal A}_a^b$ have the property that
${\tilde E}_a(\mathcal{A}^c_b)=C^c_{ad}\mathcal{A}^d_b$.

We revert to consideration of a standard basis.  We define the
component 1-forms $\omega^a$ of the tensor field $\omega$ by
$\omega=\omega^a\tilde{E}_a$. Then $\omega^a(X_i)=0$,
$\omega^a(\tilde{E}_b)=\delta^a_b$. Thus the $\omega^a$ comprise part
of the basis of 1-forms dual to the standard basis. We denote by
$\vartheta^i$ the remaining 1-forms in the dual basis.

Most of the objects of interest, such as the Lagrangian and the
corresponding Euler-Lagrange field $\Gamma$, live on the tangent
manifold of $M$, which we denote by $\tau:TM\to M$.  We recall that
there are two canonical ways of lifting a vector field, say $Z$, from
$M$ to $TM$.  The first is the complete or tangent lift, $\clift{Z}$,
whose flow consists of the tangent maps of the flow of $Z$.  The
second is the vertical lift, $\vlift{Z}$, which is tangent to the
fibres of $\tau$ and on the fibre over $m$ coincides with the constant
vector field $Z_m$.  We have $T\tau(\clift{Z})=Z$ while
$T\tau(\vlift{Z})=0$.  Moreover, $TM$ is equipped with a canonical
type $(1,1)$ tensor field called the vertical endomorphism and denoted
by $S$, which is essentially determined by the facts that $S(\clift
Z)=\vlift Z$ and $S(\vlift Z)=0$.  For more details on this material,
see e.g.\ \cite{CP,YI}.  The set $\{\clift{X_i},\clift{\tilde{E}_a},
\vlift{X_i},\vlift{\tilde{E}_a}\}$, consisting of the complete and
vertical lifts of $\{X_i,{\tilde E}_a\}$, forms a local basis of
vector fields on $TM$.

Let $\{Z_\alpha\}$ be a local basis of vector fields on $M$, and
$\{\theta^\alpha\}$ the dual basis of 1-forms.  These 1-forms define
fibre-linear functions $\vec{\theta}^a$ on $TM$, such that for any
$u\in T_mM$, $u=\vec{\theta}^a(u)Z_\alpha(m)$.  These functions are
therefore the components of velocities with respect to the specified
vector-field basis. We may use these functions as fibre coordinates.
Coordinates of this type are sometimes called quasi-velocities, and we
will use this terminology. In the case of interest we have a
standard basis $\{X_i,{\tilde E}_a\}$ and its dual
$\{\vartheta^i,\omega^a\}$; we denote the corresponding
quasi-velocities by $v^i=\vec{\vartheta}^i$, $v^a=\vec{\omega}^a$.

We will need to evaluate the actions of the vector fields
$\clift{X_i}$, $\clift{\tilde{E}_a}$, $\vlift{X_i}$ and
$\vlift{\tilde{E}_a}$ on $v^i$ and $v^a$.  Now for any vector field
$Z$ and 1-form $\theta$ on $M$,
\[
\clift{Z}(\vec{\theta})=\overrightarrow{\lie{Z}\theta}, \quad
\vlift{Z}(\vec{\theta})=\tau^*\theta(Z).
\]
Most of the required results are easy to derive from these formulae.
The only tricky calculation is that of $\clift{X_i}(\vec{\omega}^a)$,
for which we need the Lie derivative of a connection form by a
horizontal vector field.  We have
\begin{eqnarray*}
(\lie{X_i}\omega^a)(\tilde{E}_b)&=&
X_i(\delta^a_b)-\omega^a([X_i,\tilde{E}_b])=0,\\
(\lie{X_i}\omega^a)(X_j)&=&-\omega^a([X_i,X_j])=-R^a_{ij};
\end{eqnarray*}
in the first we have used the invariance of the horizontal vector
fields.  In summary, the relevant derivatives of the
quasi-velocities are
\[
\begin{array}{lllllll}
\clift{X_i}(v^j)=0,&\quad&\vlift{X_i}(v^j)=\delta^j_i,&\quad&
\clift{X_i}(v^a)=-R^a_{ij}v^j,&\quad&\vlift{X_i}(v^a)=0,\\
\clift{\tilde{E}_a}(v^i)=0,&&\vlift{\tilde{E}_a}(v^i)=0,&&
\clift{\tilde{E}_a}(v^b)=C_{ac}^bv^c,&&\vlift{\tilde{E}_a}(v^b)=\delta^b_a.
\end{array}
\]
Finally, we list some important Lie brackets of the basis vector
fields:
\begin{eqnarray*}
[\clift{{\tilde E}_a}, \clift{X_i}] = \clift{[{\tilde E}_a, X_i]}= 0,
&\quad&
[\clift{{\tilde E}_a}, \vlift{X_i}] = \vlift{[{\tilde E}_a, X_i]} = 0,\\
\mbox{}[\clift{{\tilde E}_a}, \clift{{\tilde E}_b}]
= \clift{[{\tilde E}_a, {\tilde E}_b]}
= -C^c_{ab} \clift{{\tilde E}_{c}},&\quad&
[\clift{\tilde{E}_a},\vlift{\tilde{E}_b}]=\vlift{[\tilde{E}_a,\tilde{E}_b]}
=-C_{ab}^c\vlift{\tilde{E}_c}.
\end{eqnarray*}

\section{The generalized Routh equations}

We begin by explaining, in general terms, how we will deal with the
Euler-Lagrange equations.

Consider a manifold $M$, with local coordinates $(x^\alpha)$, and
its tangent bundle $\tau:TM\to M$, with corresponding local
coordinates $(x^\alpha,u^\alpha)$. A Lagrangian $L$ is a function on
$TM$; its Euler-Lagrange equations,
\[
\frac{d}{dt}\left(\fpd{L}{u^\alpha}\right)-\fpd{L}{x^\alpha}=0,
\]
comprise a system of second-order ordinary differential equations
for the extremals; in general the second derivatives
$\ddot{x}^\alpha$ are given implicitly by these equations. We say
that $L$ is regular if its Hessian with respect to the fibre
coordinates,
\[
\spd{L}{u^\alpha}{u^\beta},
\]
considered as a symmetric matrix, is everywhere non-singular.  When
the Lagrangian is regular the Euler-Lagrange equations may be solved
explicitly for the $\ddot{x}^\alpha$, and so determine a system of
differential equations of the form
$\ddot{x}^\alpha=f^\alpha(x,\dot{x})$.  These equations can in turn
be thought of as defining a vector field $\Gamma$ on $TM$, a
second-order differential equation field, namely
\[
\Gamma=u^\alpha\vf{x^\alpha}+f^\alpha\vf{u^\alpha};
\]
we call this the Euler-Lagrange field of $L$. The Euler-Lagrange
equations may be written
\[
\Gamma\left(\fpd{L}{u^\alpha}\right)-\fpd{L}{x^\alpha}=0,
\]
and when $L$ is regular these equations, together with the
assumption that it is a second-order differential equation field,
determine $\Gamma$.

This is essentially how we will deal with the Euler-Lagrange
equations throughout:\ that is, we will assume that $L$ is regular
and we will work with the Euler-Lagrange field $\Gamma$, and with the
Euler-Lagrange equations in the form given above. However, we need to
be able to express those equations in terms of a basis of
vector fields on $M$ which is not necessarily of coordinate type. It
is easy to see that if $\{Z_\alpha\}$ is such a basis then the
equations
\[
\Gamma(\vlift{Z_\alpha}(L))-\clift{Z_\alpha}(L)=0
\]
are equivalent to the Euler-Lagrange equations. The fact that
$\Gamma$ is a second-order differential equation field means that it
takes the form
\[
\Gamma=w^\alpha\clift{Z_\alpha}+ {\Gamma}^\alpha\vlift{Z_\alpha}
\]
where the $w^\alpha$ are the quasi-velocities corresponding to the
basis $\{Z_\alpha\}$.

We now build in the assumption that $L$ has a symmetry group $G$,
which acts in such a way that $M$ is a principal bundle with $G$ as
its group, as we described above.  We will suppose that the
Lagrangian is invariant under the induced action of $G$ on $TM$.
This tangent action is defined by the collection of transformations
$\psiTM_g=T\psiM_g$ on $TM$, $g\in G$.  By construction, the
fundamental vector fields for this induced action are the complete
lifts of the fundamental vector fields of the action on $M$; the
invariance of the Lagrangian can therefore be characterized by the
property $\clift{\tilde E}_a(L)=0$.  We have shown in \cite{MCPaper}
that if $L$ is invariant, then so also is $\Gamma$, which is to say
that $[\clift{\tilde E}_a,\Gamma]=0$.

We choose a principal connection on $M$, and a basis of vector
fields $\{X_i,\tilde{E}_a\}$ adapted to it (a standard basis), as
described above.  Then the Euler-Lagrange equations for $L$ are
\begin{eqnarray*}
\Gamma(\vlift{X_i}(L))-\clift{X_i}(L)&=&0\\
\Gamma(\vlift{\tilde{E}_a}(L))-\clift{\tilde{E}_a}(L)&=&0.
\end{eqnarray*}
But by assumption $\clift{\tilde{E}_a}(L)=0$:\ it follows
immediately that $\Gamma(\vlift{\tilde{E}_a}(L))=0$.  So the
functions $\vlift{\tilde{E}_a}(L)$ are first integrals, which
clearly generalize the momenta conjugate to ignorable coordinates in
the classical Routhian picture.  We write $p_a$ for
$\vlift{\tilde{E}_a}(L)$.  The Euler-Lagrange field is tangent to
any submanifold $p_a=\mu_a=\mbox{constant}$,
$a=1,2,\ldots,\dim(\g)$, that is, any level set of momentum.  By a
well-known argument (see e.g.\ \cite{MRbook}), we may regard
$(x,v)\mapsto(p_a(x,v))$ as a map from $TM$ to $\g^*$, the dual of
the Lie algebra $\g$, and this map is equivariant between the given
action of $G$ on $TM$ and the coadjoint action of $G$ on $\g^*$ (the
coadjoint action is defined as $\langle\xi,\ad^*_g\mu\rangle
=\langle\ad_g\xi,\mu\rangle$).  We have
\[
\clift{\tilde{E}_a}(p_b)=\clift{\tilde{E}_a}\vlift{\tilde{E}_b}(L)
=[\clift{\tilde{E}_a},\vlift{\tilde{E}_b}](L)
=-C_{ab}^c\vlift{\tilde{E}_c}(L)=-C_{ab}^cp_c,
\]
which expresses this result in our formalism.

We will also need a less coordinate-dependent version of the
Hessian. In fact the Hessian of $L$ at $w\in TM$ is the symmetric
bilinear form $g$ on $T_mM$, $m=\tau(w)$, given by
$g(u,v)=\vlift{u}\vlift{v}(L)$, where the vertical lifts are to $w$.
We can equally well regard $g$ as a bilinear form on the vertical
subspace of $T_wTM$, by identifying $u$ and $v$ with their vertical
lifts. The components of the Hessian
$g$ with respect to our standard basis will be denoted as follows:
\[
g(\tilde{E}_a,\tilde{E}_b)=g_{ab},\quad g(X_i,X_j)=g_{ij},\quad
g(X_i,\tilde{E}_a)=g_{ia}=g_{ai}=g(\tilde{E}_a,X_i).
\]
We also have $g_{ab}=\vlift{\tilde{E}_a}(p_b)$,
$g_{ia}=\vlift{X_i}(p_a)$. In general these components are functions
on $TM$, not on $M$, and the Hessian should be regarded as a tensor
field along the tangent bundle projection $\tau: TM \to M$. We will
assume throughout that $L$ is regular, which means that $g$ as a
whole is non-singular.  Then $\Gamma$ is uniquely determined as a
second-order differential equation field on $TM$.

We now turn to the consideration of Routh's procedure.  We call the
function $\R$ on $TM$ given by
\[
\R=L-v^ap_a
\]
the {\em Routhian}. It generalizes in an obvious way the classical
Routhian corresponding to ignorable coordinates. The Routhian is
invariant:
\[
\clift{\tilde{E}_b}(v^ap_a)=
C^a_{bc}v^cp_a-v^aC^c_{ba}p_c=0,
\]
whence the result.

We now consider the Euler-Lagrange equations
$\Gamma(\vlift{X_i}(L))-\clift{X_i}(L)=0$.  We wish to write these
equations in terms of the restriction of the Routhian to a level set
of momentum, say $p_a=\mu_a$, which we denote by $N_\mu$.  To do so,
we need to work in terms of vector fields related to $\clift{X_i}$,
$\vlift{X_i}$ and $\clift{\tilde{E}_a}$ which are tangent to $N_\mu$
(in general there is no reason to suppose that these vector fields
themselves have this property, of course).  To define the new vector
fields we will assume that the Lagrangian has an additional
regularity property: we will assume that $(g_{ab})$ is non-singular.
(Note that if the Hessian is everywhere positive-definite then
$(g_{ab})$ is automatically non-singular.) Then there are
coefficients $A^b_i$, $B^b_i$ and $C^b_a$, uniquely defined, such
that
\begin{eqnarray*}
(\clift{X_i}+A^b_i\vlift{\tilde{E}_b})(p_a)&=&
\clift{X_i}(p_a)+A^b_ig_{ab}=0\\
(\vlift{X_i}+B^b_i\vlift{\tilde{E}_b})(p_a)&=&
\vlift{X_i}(p_a)+B^b_ig_{ab}=0\\
(\clift{\tilde{E}_a}+C^b_a\vlift{\tilde{E}_b})(p_c)&=&
\clift{\tilde{E}_a}(p_c)+C^b_ag_{bc}=0.
\end{eqnarray*}
The vector fields $\clift{\bar{X}_i}$, $\vlift{\bar{X}_i}$ and
$\clift{\bar{E}_a}$ given by
\begin{eqnarray*}
\clift{\bar{X}_i}&=&\clift{X_i}+A^a_i\vlift{\tilde{E}_a}\\
\vlift{\bar{X}_i}&=&\vlift{X_i}+B^a_i\vlift{\tilde{E}_a}\\
\clift{\bar{E}_a}&=&\clift{\tilde{E}_a}+C^b_a\vlift{\tilde{E}_b}
\end{eqnarray*}
are tangent to each level set $N_\mu$.  (The notation is not meant
to imply that the barred vector fields are actually complete or
vertical lifts.)  We will need to know the coefficients explicitly
only in the case of $B^a_i$ and $C^b_a$:\ in fact
\[
B^a_i=-g^{ab}g_{ib} \qquad\mbox{and} \qquad
C^b_a=g^{bc}C^d_{ac}p_d.
\]

This is all carried out under the assumption that $(g_{ab})$ is
non-singular.  One has to make such an assumption in the classical
case in order to be able to solve the equations $\partial
L/\partial\dot{\theta}^a=\mu_a$ for the $\dot{\theta}^a$.  In the
general case the non-singularity of $(g_{ab})$ is the condition for
the level set $N_\mu$ to be regular, i.e.\ to define a submanifold
of $TM$ of codimension $\dim(\g)$.  The vector fields
$\vlift{\tilde{E}_a}$ are transverse to all regular level sets, and
the barred vector fields span the level sets.  Thus on any regular
level set the bracket of any two of the barred vector fields is a
linear combination of vector fields of the same form.  We want in
particular to observe that this implies that
$[\clift{\bar{E}_a},\vlift{\bar{X}_i}]=0$.  It is not difficult to
see, using the known facts about the brackets of the unbarred vector
fields, that this bracket is of the form $P^a\vlift{\tilde{E}_a}$;
this must satisfy $P^a\vlift{\tilde{E}_a}(p_b)=P^ag_{ab}=0$, whence
by the regularity assumption $P^a=0$.  In fact by similar arguments
the brackets of the barred vector fields just reproduce those of
their unbarred counterparts, except that $[\clift{\bar
X}_i,\vlift{\bar X}_j]=0$.  In particular,
$[\clift{\bar{E}_a},\clift{\bar{E}_b}]=-C^c_{ab}\clift{\bar{E}_c}$.
The $\clift{\bar{E}_a}$ therefore form an anti-representation of
$\g$, acting on the level set $N_\mu$ (just as the $\tilde{E}_a$ do
on $M$).

We return to the expression of the Euler-Lagrange equations in terms
of the Routhian. We will need to evaluate the actions of
$\clift{\bar{X}_i}$ and $\vlift{\bar{X}_i}$ on $v^a$.
Using the formulae in Section~2 we find that
\begin{eqnarray*}
\clift{\bar{X}_i}(v^a)&=&
(\clift{X_i}+A^b_i\vlift{\tilde{E}_b})(v^a)=
-R^a_{ij}v^j+A^a_i\\
\vlift{\bar{X}_i}(v^a)&=&
(\vlift{X_i}+B^b_i\vlift{\tilde{E}_b})(v^a)=B^a_i.
\end{eqnarray*}
We now set things up so that we can restrict to the submanifold
$N_\mu$ easily. We have
\begin{eqnarray*}
\clift{X_i}(L)&=&\clift{\bar{X}_i}(L)-A^a_i\vlift{\tilde{E}_a}(L)\\
&=&\clift{\bar{X}_i}(L-v^ap_a)+(-R^a_{ij}v^j+A^a_i)p_a
+v^a\clift{\bar{X}_i}(p_a)-A^a_ip_a\\
&=&\clift{\bar{X}_i}(\R)-p_aR^a_{ij}v^j;\\
\vlift{X_i}(L)&=&\vlift{\bar{X}_i}(L)-B^a_i\vlift{\tilde{E}_a}(L)\\
&=&\vlift{\bar{X}_i}(L-v^ap_a)+B^a_ip_a
+v^a\vlift{\bar{X}_i}(p_a)-B^a_ip_a\\
&=&\vlift{\bar{X}_i}(\R).
\end{eqnarray*}
But $\Gamma(\vlift{X_i}(L))-\clift{X_i}(L)=0$, and $\Gamma$ is tangent
to the submanifold $N_\mu$; thus if we denote
by $\R^\mu$ the restriction of the Routhian to the submanifold (where
it becomes $L-v^a\mu_a$) we have
\[
\Gamma(\vlift{\bar{X}_i}(\R^\mu))-\clift{\bar{X}_i}(\R^\mu)=-\mu_aR^a_{ij}v^j.
\]
On the other hand, if $\Gamma$ is a second-order differential equation
field such that $\Gamma(\vlift{\tilde{E}_a}(L))=0$ and the
above equation holds for all $\mu_a$ then $\Gamma$ satisfies the
Euler-Lagrange equations for the invariant Lagrangian $L$.

We will refer to these equations as the generalized Routh
equations.

Neither $\R^\mu$ nor $\Gamma$ is $\clift{\bar{E}_a}$-invariant.  They
will however be invariant under those vector fields
$\clift{\tilde{\xi}}$, $\xi\in\g$, which happen to be tangent to the
level set $N_\mu$.  These are the vector fields for which
$\xi^a\clift{\tilde{E}_a}=\xi^a\clift{\bar{E}_a}$, or
$\xi^aC_{ab}^c\mu_c=0$. We will return to this issue in later
sections.

Note that since $\Gamma$ satisfies $\Gamma(p_a)=0$ it may be
expressed in the form
\[
\Gamma=v^i\clift{\bar{X}_i}+\Gamma^i\vlift{\bar{X}_i}+v^a\clift{\bar{E}_a}.
\]
If the matrix-valued function $\vlift{\bar X}_i \vlift{\bar X}_j(\R)$
is non-singular, the reduced Euler-Lagrange equations above will
determine the coefficients $\Gamma^i$.  We show now that this is
always the case, under the assumptions made earlier.

Recall that $\vlift{\bar{X}_i}=\vlift{X_i}+B^a_i\vlift{\tilde{E}_a}$
is determined by the condition that $\vlift{\bar{X}_i}(p_a)=0$, and
that therefore $B^a_i=-g^{ab}g_{ib}$. We may regard
$X_i+B^a_i\tilde{E}_a$ as a vector field along the tangent bundle
projection, and $\vlift{\bar{X}_i}$ really is the vertical lift of
this vector field; we will accordingly denote it by $\bar{X}_i$.
Then
\[
g(\bar{X}_i,\tilde{E}_a)=g(X_i,\tilde{E}_a)+B^b_ig(\tilde{E}_b,\tilde{E}_a)=
g_{ia}+B^b_ig_{ab}=0;
\]
thus the $\bar{X}_i$ span the orthogonal complement to the space
spanned by the $\tilde{E}_a$ with respect to the Hessian of $L$.
That is to say, the tangent space to a regular level set of momentum
at any point $u\in TM$ intersects the tangent space to the fibre of
$TM\to M$ at $u$ in the subspace orthogonal with respect to $g_u$ to
the span of the $\vlift{\tilde{E}_a}$. Moreover,
\begin{eqnarray*}
g(\bar{X}_i,\bar{X}_j)&=&
g_{ij}+B^a_ig_{aj}+B^a_jg_{ia}+B^a_iB^b_jg_{ab}\\
&=&g_{ij}-2g^{ab}g_{ia}g_{jb}+g^{ac}g^{bd}g_{ic}g_{jd}g_{ab}\\
&=&g_{ij}-g^{ab}g_{ia}g_{jb}.
\end{eqnarray*}
So this is the expression for the restriction of the Hessian of $L$
to the subspace orthogonal to that spanned by the $\tilde{E}_a$.

Now recall that $\vlift{\bar{X}_i}(\R)=\vlift{X_i}(L)$. Thus
\[
\vlift{\bar{X}_i}\vlift{\bar{X}_j}(\R)=
(\vlift{X_i}-g^{ab}g_{ib}\vlift{\tilde{E}_a})\vlift{X_j}(L)
=g_{ij}-g^{ab}g_{ib}g_{aj}=g(\bar{X}_i,\bar{X}_j).
\]
That is, the `Hessian' of $\R$ (i.e.\
$\vlift{\bar{X}_i}\vlift{\bar{X}_j}(\R)$) is just the restriction of
the Hessian of $L$ to the subspace orthogonal to that spanned by the
$\tilde{E}_a$.  It follows that the bilinear form with components
$\bar{g}_{ij}=\vlift{\bar{X}_i}\vlift{\bar{X}_j}(\R)$ must be
non-singular. For suppose that there is some vector $w^j$ such that
$\bar{g}_{ij}w^j=0$; then $g(\bar{X}_i,w^j\bar{X}_j)=0$ by
assumption, and $g(\tilde{E}_a,w^j\bar{X}_j)=0$ by orthogonality ---
but then $w^j\bar{X}_j=0$ since $g$ is assumed to be non-singular.

The sense in which the generalized Routh equations are `reduced'
Euler-Lagrange equations is that (in principle at least) we can reduce
the number of variables by using the equations $p_a=\mu_a$ to
eliminate the quasi-velocities $v^a$.  However, these variables appear
explicitly in the expression for $\Gamma$, so it may be considered
desirable to rearrange the generalized Routh equations so that they no
longer appear.  This can be done by changing the basis of vector
fields on the level set of momentum, as follows.  The change is
suggested by the fact that, notation notwithstanding,
$S(\clift{\bar{X}_i})\neq\vlift{\bar{X}_i}$ (where $S$ is the vertical
endomorphism).  Let us, however, set
\[
\clift{\hat{X}_i}=\clift{\bar{X}_i}+B^a_i \clift{\bar{E}_a}:
\]
then since $S$ vanishes on vertical lifts,
\[
S(\clift{\hat{X}_i})=
S(\clift{X_i}+B^a_i\clift{\tilde{E}_a})=
\vlift{X_i}+B^a_i\vlift{\tilde{E}_a}=\vlift{\bar{X}_i}.
\]
We write
\[
\Gamma_0=v^i\clift{\hat{X}_i}+\Gamma^i\vlift{\bar{X}_i},
\]
so that
\[
\Gamma=\Gamma_0+(v^iB^a_i+v^a)\clift{\bar{E}_a}.
\]
We will examine the contribution of the term involving
$\clift{\bar{E}_a}$ in $\Gamma$ to the generalized
Routh equations. First we determine
$\clift{\bar{E}_a}(\R)$.  Since $\clift{\bar{E}_a}(p_b)=0$,
\[
\clift{\bar{E}_a}(\R)=\clift{\bar{E}_a}(L-v^bp_b)
=C_a^bp_b-\clift{\bar{E}_a}(v^b)p_b
=C_a^bp_b-C_{ac}^bp_bv^c-C_a^c\delta_c^bp_b=-C_{ac}^bp_bv^c.
\]
It follows that
\[
\clift{\bar{E}_a}(\vlift{\bar{X}_i}(\R))=
\vlift{\bar{X}_i}(\clift{\bar{E}_a}(\R))=
-\vlift{\bar{X}_i}(C_{ac}^bp_bv^c) = -C_{ac}^bp_b B^c_i.
\]
So setting $\Gamma=\Gamma_0 + (v^iB^a_i+v^a)\clift{\bar{E}_a}$ we have
\begin{eqnarray*}
\Gamma(\vlift{\bar{X}_i}(\R^\mu))-\clift{\bar{X}_i}(\R^\mu)&=&
\Gamma_0(\vlift{\bar{X}_i}(\R^\mu))
+(v^jB^a_j+v^a)\clift{\bar{E}_a}(\vlift{\bar{X}_i}(\R^\mu))
-\clift{\hat{X}_i}(\R^\mu)+B^a_i \clift{\bar{E}_a}(\R^\mu)\\
&=&
\Gamma_0(\vlift{\bar{X}_i}(\R^\mu))-\clift{\hat{X}_i}(\R^\mu)
-(v^jB^a_j+v^a)C_{ac}^b\mu_b B^c_i-B^a_iC_{ac}^b\mu_bv^c\\
&=&
\Gamma_0(\vlift{\bar{X}_i}(\R^\mu))-\clift{\hat{X}_i}(\R^\mu)
+v^jB^a_jB^c_iC_{ac}^b\mu_b,
\end{eqnarray*}
and the generalized Routh equations become
\[
\Gamma_0(\vlift{\bar{X}_i}(\R^\mu))-\clift{\hat{X}_i}(\R^\mu)=
-\mu_a(R^a_{ij}+B^b_iB^c_jC_{bc}^a)v^j.
\]
We may say that among the vector fields tangent to a level set of
momentum it is $\clift{\bar{X}_i}+B^a_i \clift{\bar{E}_a}$, not
$\clift{\bar{X}_i}$, that really plays the role of the complete lift
of $\bar{X}_i$. Be aware, however, that unless the symmetry group is
Abelian, $\Gamma_0$ cannot be identified with a vector field on
$T(M/G)$. We will return to this matter at the end of
Section~\ref{conn}.

To end this section we give a coordinate expression for the
generalized Routh equations in their original form. For this purpose
we take coordinates $(x^i)$ on $M/G$, and coordinates
$(x^i,\theta^a)$ on $M$ such that the $\theta^a$ are fibre
coordinates; then $(x^i,\theta^a,v^i)$ are coordinates on $N_\mu$,
which is to say that $N_\mu$ can be locally identified with
$M\times_{M/G}T(M/G)$.  We may write
\[
X_i=\vf{x^i}-\Lambda^a_i\vf{\theta^a},\quad
\tilde{E}_a=K^b_a\vf{\theta^b}
\]
for suitable functions $\Lambda^a_i$ and $K^b_a$ on $M$.  (We should
note that the $K^b_a$ are components of a non-singular matrix at each
point; moreover, the invariance property of the $X_i$ can be expressed
in terms of the coefficients $\Lambda^a_i$ and $K^b_a$; but we will
not actually need either of these facts here.)  From the formulae for
the action of complete and vertical lifts on quasi-velocities given at
the end of Section 2 we see that
\[
\clift{\bar{X}_i}(v^i)=\clift{\bar{E}_a}(v^i)=0,\quad
\vlift{\bar{X}_i}(v^j)=\delta^j_i.
\]
Thus in terms of $x^i$, $\theta^a$ and $v^i$ we can write
\[
\clift{\bar{X}_i}=\vf{x^i}-\Lambda^a_i\vf{\theta^a},\quad
\vlift{\bar{X}_i}=\vf{v^i},\quad
\clift{\bar{E}_a}=K^b_a\vf{\theta^b}.
\]
It is necessary to be a little careful:\ the coordinate vector field
expressions are ambiguous, since they can refer either to coordinates
on $TM$ or on $N_\mu$.  We emphasise that it is the latter
interpretation that is intended here. In view of the possibilities of
confusion it will be useful to have an explicit notation for the
injection $\J\to TM$:\ we denote it by $\iota$. The non-singularity of $(g_{ab})$ ensures
that, at least locally, we can rewrite the relation $p_a=\mu_a$ for
the injection $\iota:\J\to TM$ in the form $v^a =
\iota^a(x^i,\theta^a,v^i)$, for certain functions $\iota^a$ of the
specified variables.

The restriction of the Euler-Lagrange field $\Gamma$ to
$\J$ is
\begin{eqnarray*}
\Gamma &=&
\iota^a\clift{\bar{E}_a}+ v^i\clift{\bar{X}_i}+
(\Gamma^i\circ\iota) \vlift{\bar{X}_i}\\
&=&\iota^bK^a_b\vf{\theta^a}+v^i\left(\vf{x^i}-\Lambda^a_i\vf{\theta^a}\right)
+(\Gamma^i\circ\iota)\vf{v^i}\\
&=&\left(\iota^bK^a_b-v^i\Lambda^a_i\right)\vf{\theta^a}
+v^i\vf{x^i}+(\Gamma^i\circ\iota)\vf{v^i};
\end{eqnarray*}
the equations for its integral curves are
\[
\left\{\begin{array}{lll}
{\dot x}^i &=& v^i, \\
{\dot v}^i  &=&\Gamma^i(x,\theta,v),\\
{\dot \theta}^a &=& \iota^b(x,\theta,v) K^a_b(x,\theta) -
{v}^i{\Lambda}^a_i(x,\theta).\end{array} \right.
\]
These can be considered as a coupled set of first- and  second-order
equations,
\[
\left\{\begin{array}{lll}
{\ddot x}^i &=& \Gamma^i(x,\theta,\dot{x}),\\
{\dot \theta}^a &=& \iota^b(x,\theta,\dot{x}) K^a_b(x,\theta) -
\dot{x}^i{\Lambda}^a_i(x,\theta).\end{array} \right.
\]
With regard to the second of these equations, we point out that the
expression for the velocity variables $\dot{\theta}^a$ in terms of the
quasi-velocities $v^i$ and $v^a$ is just
$\dot{\theta}^a=v^bK^a_b-v^i\Lambda^a_i$.  What turns these
identities into genuine differential equations is, in particular,
substitution for the $v^a$ in terms of the other variables via the
functions $\iota^a$ --- or in other words, restriction to $N_\mu$.

The functions $\Gamma^i$ may be determined from the generalized Routh
equations. These may be expressed as
\[
\frac{d}{dt}\left(\fpd{\R^\mu}{v^i}\right)-\fpd{\R^\mu}{x^i}=
-\mu_aR^a_{ij}-\Lambda^a_i\fpd{\R^\mu}{\theta^a}.
\]
In the light of the earlier remarks about the interpretation of
coordinate vector fields, we point out that substitution for $v^a$ in
terms of the other variables in this equation must be carried out
before the partial derivatives are calculated.

\section{The reconstruction method}

We have seen in Section 3 that Routh's technique consists in
restricting the Euler-Lagrange equations to a level set of momentum
$N_\mu$.  This procedure takes partial, but not necessarily complete,
account of the action of the symmetry group $G$.  To make further
progress we must examine the residual action of $G$ on $N_\mu$.

As we mentioned before, the momentum map is equivariant between the
induced action of $G$ on $TM$ and the coadjoint action of $G$ on
$\g^*$.  The submanifold $N_\mu$ is therefore invariant under the
isotropy group $G_\mu= \{g\in G\mid \ad^*_g\mu=\mu\}$ of $\mu$.  The
algebra $\g_\mu$ of $G_\mu$ consists of those $\xi\in\g$ such that
$\xi^bC^c_{ab}\mu_c=0$; this is the necessary and sufficient condition
for $\clift{\tilde{\xi}}$ to be tangent to $N_\mu$.

Note that any geometric object we know to be $G$-invariant is
automatically $G_\mu$-invariant.

The manifold $N_\mu$ is a principal fibre bundle with group $G_\mu$;
we will denote its base by $N_\mu/G_\mu$.  The restriction of the
Euler-Lagrange field $\Gamma$ to $N_\mu$ is $G_\mu$-invariant, and as
a consequence it projects onto a vector field $\check\Gamma$ on
$N_\mu/G_\mu$.

The task now is to examine the relationship between $\check\Gamma$
and $\Gamma$. There are two aspects:\ the formulation of the
differential equations represented by $\check\Gamma$; and the
reconstruction of integral curves of $\Gamma$ from integral curves of
$\check\Gamma$ (supposing that we have solved those equations).

Our methods of attack on these problems will be based on those we
developed in our papers~\cite{Paper1,MCPaper} and are similar to
(but different from) the ones that were adopted in e.g.\ \cite{MMR}.
These in turn were based on the following well-known method for
reconstructing integral curves of an invariant vector field from
reduced data. Let $\pi:N\to B$ be a principal fibre bundle with
group $G$.  Any invariant vector field $\Gamma$ on $N$ defines a
$\pi$-related reduced vector field $\check \Gamma$ on $B$: due to
the invariance of $\Gamma$, the relation $T\pi \big(\Gamma (n)\big)=
\check\Gamma \big(\pi(n)\big)$ is independent of the choice of $n\in
N$ within the equivalence class of $\pi(n)\in B$. Given a principal
connection $\Omega$, an integral curve $v(t)$ of $\Gamma$ can be
reconstructed from an integral curve ${\check v}(t)$ of ${\check
\Gamma}$ as follows.  Let $\hlift{\check{v}}(t)$ be a horizontal
lift of ${\check v}(t)$ with respect to $\Omega$ (that is, a curve
in $N$ over $\check v$ such that $\Omega(\dot{\hlift{\check v}})=0$)
and let $g(t)$ be the solution in $G$ of the equation
\[
\widetilde{\vartheta(\dot g(t))} = \Omega( \Gamma(\hlift{\check v}(t)))
\]
where $\vartheta$ is the Maurer-Cartan form of $G$.  (We use here
the fact that given any curve $\xi(t)$ in $\g$, the Lie algebra of
$G$, there is a unique curve $g(t)$ in $G$ which satisfies
$\vartheta(\dot{g}(t))=\xi(t)$ and $g(0)=e$; $g(t)$ is sometimes
called the development of $\xi(t)$ into $G$, see for example
\cite{Sharpe}.) Then
$v(t)=\psiN_{g(t)}\hlift{\check{v}}(t)$ is an integral curve of
$\Gamma$.

In the following sections we define two principal connections on
$\J$, we determine $\check\Gamma$ and we identify for both
connections the vertical part of $\Gamma$, necessary for the
reconstruction method above.

\section{Two principal connections on a level set of
momentum}\label{conn}

A principal connection $\Omega$ on $\J\to\J/G_\mu$ is by definition
a left splitting of the short exact sequence
\[
0\to\J\times \g_\mu \to T\J \to \J \times_{N_\mu/G_\mu}
T(\J/G_\mu)\to 0;
\]
all spaces in the above sequence should be interpreted as bundles
over $\J$. We think of $\Omega$ as a type $(1,1)$ tensor field on
$\J$ which is pointwise a projection operator with image the tangent
space to the fibre, and which is invariant under $G_\mu$.

The first connection we define uses the Hessian of $L$ to determine
its horizontal distribution, and is therefore analogous to the
mechanical connection of a simple system; we denote it by
$\mech\Omega$.

Recall that we interpret the Hessian $g$ of $L$ as a tensor field
along $\tau$. In particular, its components with respect to the standard
basis $\{X_i,\tilde{E}_a\}$ are functions on $TM$. We will say that
a vector field $W$ on $N_\mu$ is horizontal for $\mech\Omega$ if
\[
g({\tilde\xi},\tau_* W) =0, \qquad \forall \xi\in\g_\mu,
\]
where $\tau_*W$ is the projection of a vector field $W$ on $TM$ to a
vector field along $\tau:TM\to M$.  The definition makes sense only if
we assume that the restriction of $g$ to $N_\mu\times \g_\mu$ is
non-singular, as we do from now on.

In \cite{MCPaper} we have shown that if the Lagrangian is
invariant then so is $g$, in the sense that
\[
\lie{\tilde\xi}g =0,\qquad \forall \xi\in\g.
\]
Here, for a vector field $Z$ on $M$,
$\lie{Z}$ stands for an operator acting on tensor fields along
$\tau$ that has all the properties of a Lie derivative operator, and
in particular, when applied to a function $f$ on $TM$ and a vector
field $X$ along $\tau$ gives
\[
\lie{Z}f = \clift{Z}(f),\qquad
\lie{Z}X=
\left(Z^\beta\fpd{X^\alpha}{x^\beta}+
\fpd{Z^\beta}{x^\gamma}{u}^\gamma\fpd{X^\alpha}{u^\beta}
-X^\beta\fpd{Z^\alpha}{x^\beta}\right)\vf{x^\alpha}
\]
where
\[
Z=Z^\alpha(x)\fpd{}{x^\alpha},\qquad
X=X^\alpha(x,u)\fpd{}{x^\alpha}.
\]
Note that if $X$ is a basic vector field along $\tau$ (i.e.\ a vector
field on $M$), then $\lie{Z}X = [Z,X]$.  Furthermore, for any vector
field $W$ on $TM$ we have
\[
\lie{Z}(\tau_*W)=\tau_*[\clift Z,W].
\]

To show that the connection is principal we need only to show that
if $W$ is horizontal so also is $[\clift{\tilde\xi},W]$ for all
$\xi\in\g_\mu$. But for all $\xi,\eta\in\g_\mu$,
\[
g(\tau_*[\clift{\tilde{\xi}},W],\tilde{\eta})=
g(\lie{{\tilde{\xi}}}(\tau_*W),\tilde{\eta})
=-g(\tau_*W,\lie{{\tilde{\xi}}}\tilde{\eta})=
-g(\tau_*W,[\tilde{\xi},\tilde{\eta}])=
g(\tau_*W,\widetilde{[\xi,\eta]})=0,
\]
using the properties of the generalized Lie derivative and the
invariance of $g$.

As was mentioned before, $\{\clift{\bar X}_i,\clift{\bar
E}_a,\vlift{\bar X}_i\}$ is a basis of vector fields on $\J$.
Suppose now that the basis $\{E_a\} = \{E_A, E_\alpha\}$ of $\g$ is
chosen so that $\{E_A\}$ is a basis of $\g_\mu$.  Then
$C_{Ab}^c\mu_c=0$, and on $N_\mu$ we get for the corresponding
fundamental vector fields
\[
\clift{\bar E}_A = \clift{\tilde E}_A + g^{bc} C_{Ac}^d \mu_d
\vlift{\tilde E}_b = \clift{\tilde E}_A.
\]
All $\clift{\tilde E}_A$ are therefore tangent to $\J$, as required.
These vector fields span exactly the vertical space of $\J\to
\J/G_\mu$ which we have identified with $\J\times \g_\mu$. Vector
fields of this form are infinitesimal generators of the
$G_\mu$-action on $\J$.

If $(G^{AB})$ is the inverse of the matrix $(g_{AB})$ (and not the
$(A,B)$-component of $(g^{ab})$), then the vector fields
\begin{eqnarray*}\Hb_\alpha&=& \clift{\bar E}_\alpha -
G^{AB}g_{A\alpha}\clift{\tilde E}_B = \clift{\bar E}_\alpha -
\Upsilon^B_\alpha \clift{\tilde E}_B
\\\hlift{\bar X}_i&=&\clift{\bar X}_i - G^{AB}g_{Ai}\clift{\tilde
E}_B = \clift{\bar X}_i - \Upsilon^B_i \clift{\tilde E}_B,
\end{eqnarray*}
together with $\vlift{\bar X}_i$, are horizontal.  (As was the case
with the notations $\clift{{\bar E}_a}$ etc., the notation for the
horizontal fields is not meant to imply that $\Hb_\alpha$ etc.\ are
actually horizontal lifts.)  The action of $\mech\Omega$ is simply
\[
\mech\Omega(\clift{\tilde E}_A) = \clift{\tilde E}_A,\qquad
\mech\Omega(\Hb_\alpha) = 0,\qquad
\mech\Omega(\hlift{\bar X}_i)=0,\qquad
\mech\Omega(\vlift{\bar X}_i) =0,
\]
and since the arguments form a basis of vector fields on $N_\mu$
these equations specify $\mech\Omega$ explicitly. We will call
$\mech\Omega$ the mechanical connection on $N_\mu$.

The vector fields $\clift{\hat X}_i=\clift{\bar X}_i + B^a_i
\clift{\bar E}_a$ introduced earlier are also horizontal; they can be
expressed as $\clift{\hat X}_i = \hlift{\bar X}_i - g^{\alpha
b}g_{bi}\hlift{\bar E}_\alpha$.

The vector fields $\clift{\bar X}_i$ are not horizontal with respect
to $\mech\Omega$.  However, it is possible to identify a second
principal connection $\Omega^{\J}$ on $N_\mu$ for which these vector
fields are horizontal.  We will identify $\Omega^{\J}$ in two steps.

It seems natural to split the basis $\{\clift{\bar X}_i, \clift{\bar
E}_a, \vlift{\bar X}_i\}$ into a `vertical' part $\{\clift{\bar
E}_a\}$ and a `horizontal' part $\{\clift{\bar X}_i, \vlift{\bar
X}_i\}$.  To see that it does indeed make sense to do so it is
sufficient to observe that the distributions spanned by
$\{\clift{\bar E}_a\}$ and $\{\clift{\bar X}_i,\vlift{\bar X}_i\}$,
respectively, are unchanged when the bases $\{E_a\}$ of $\g$ and
$\{X_i\}$ of $\omega$-horizontal vector fields on $M$ are replaced
by different ones.  Under a change of basis for $\g$ the
$\{\clift{\bar E}_a\}$ are simply replaced by constant linear
combinations of themselves, so their span is clearly unchanged.  On
the other hand, if we set $Y_i=A_i^jX_j$ then
$\vlift{Y_i}=A_i^j\vlift{X_j}$ and
$\clift{Y_i}=A_i^j\clift{X_j}+\dot{A}_i^j\vlift{X_j}$ (where
$\dot{A}_i^j$ is the total derivative of $A_i^j$, not that it
matters), so the distributions spanned by $\{\clift{\bar
X}_i,\vlift{\bar X}_i\}$ and $\{\clift{\bar Y}_i,\vlift{\bar Y}_i\}$
are the same.

So we can indeed characterize a connection in this way, but it is not
a connection on $N_\mu\to N_\mu/G_\mu$.  In fact this construction
defines a connection on the bundle with projection $\J
\to T(M/G)$ (the restriction of $T\pM: TM\to T(M/G)$ to $\J$), i.e.\ a
splitting of the short exact sequence
\[
0\to \J \times \g \to T\J \to \J \times_{T(M/G)} T(T(M/G))\to 0.
\]
(Recall that the vector fields $\clift{\bar{E}_a}$, which span the
vertical space of the projection $\J\to T(M/G)$, form an
anti-representation of $\g$ acting on the level set $N_\mu$.)  The
construction just described is a version of the so-called vertical
lift of a connection on a principal bundle (here $\omega$) to its
tangent bundle (this is described more fully in \cite{Paper1});
accordingly we denote the corresponding type $(1,1)$ tensor field by
$\vlift\Omega$, and we have
\[
\vlift\Omega(\clift{\bar E}_a) = \clift{\bar E}_a,\qquad
\vlift\Omega( \clift{\bar X}_i)= 0,\qquad \vlift\Omega( \vlift{\bar
X}_i) =0.
\]
Evidently $(\vlift\Omega)^2=\vlift\Omega$. We  show now that
${\mathcal L}_{\clift{\tilde E}_A}\vlift\Omega=0$ for all $A$.
Firstly, note that
\[
[\clift{\tilde E}_A,\clift{\bar E}_a]=[\clift{\bar E}_A,\clift{\bar E}_a]
=-C_{Aa}^b\clift{\bar E}_b,
\]
so that
\[
({\mathcal L}_{\clift{\tilde E}_A}\vlift\Omega) (\clift{\bar E}_a)=
[\clift{\tilde E}_A, \vlift\Omega (\clift{\bar E}_a)]
-\vlift\Omega[\clift{\tilde E}_A,\clift{\bar E}_a]=
[\clift{\tilde E}_A,\clift{\bar E}_a]+C_{Aa}^b\clift{\bar E}_b = 0.
\]
Moreover, since $[\clift{\tilde E}_A,\clift{\bar X}_i]=0$,
\[
({\mathcal L}_{\clift{\tilde E}_A}\vlift\Omega) (\clift{\bar X}_i)=
[\clift{\tilde E}_A, \vlift\Omega (\clift{\bar X}_i)]-\vlift\Omega
[\clift{\tilde E}_A,\clift{\bar X}_i]= 0,
\]
and similarly for $\vlift{\bar X}_i$.

The relation between $\vlift{\Omega}$ and $\omega$ may be described
more easily if we momentarily break our convention by specifying
connections by their forms rather than by the tensors corresponding
to their splittings:\ it is easily checked that
\[
\vlift{\Omega}(Z_v)=\omega(T(\tau\circ\iota)Z_v),\qquad Z_v\in TN_\mu.
\]
This equation has to be read as one between elements of $\g$, obtained
by identifying the vertical subbundle of $TN_\mu$ with $N_\mu\times\g$
and the vertical subbundle of $TN$ with $M\times\g$, or if you will by
projection onto $\g$.

The vertical space $\J\times \g_\mu$ of the connection $\Omega^{\J}$
we are looking for is only a subbundle of the vertical space
$\J\times \g$ of the connection $\vlift\Omega$. So in a second step
we need to identify a connection for the following sequence of
trivial vector bundles:
\[
0\to \J\times \g_\mu \to \J\times \g \to \J \times \g/\g_\mu \to 0.
\]
For this connection we can simply take the restriction of the
mechanical connection $\mech\Omega$ defined earlier to the
submanifold $\J\times \g$.  The connection $\Omega^{\J}$ is then
simply $\mech\Omega\circ\vlift\Omega$ (see the diagram below).

\setlength{\unitlength}{1cm}

\begin{picture}(14.29,6.2)(-3,4.2)
\put(0.9,9.9){\vector(1,0){3.3}} \put(5.7,9.9){\vector(1,0){3.1}}

\put(0.9,7.3){\vector(1,0){3.5}} \put(5.5,7.3){\vector(1,0){2.5}}

\put(0.5,4.7){\vector(1,0){2.4}} \put(7,4.7){\vector(1,0){0.8}}

\put(9.9,9.5){\vector(0,-1){1.7}} \put(5,6.8){\vector(0,-1){1.7}}
\put(5,9.5){\vector(0,-1){1.7}} \put(0.1,6.8){\vector(0,-1){1.7}}
\put(0.1,9.5){\vector(0,-1){1.7}} \put(9.9,6.8){\vector(0,-1){1.7}}

\put(0.1,9.9){\mybox{$\J \times \g_\mu$}}
\put(5,9.9){\mybox{$\J\times \g$}} \put(9.9,9.9){\mybox{$\J \times
\g/\g_\mu$}}

\put(0.1,7.3){\mybox{$\J \times \g_\mu$}} \put(5,7.3){\mybox{$T\J$}}
\put(9.9,7.3){\mybox{$\J\times_{\J/G_\mu} T(\J/G_\mu)$}}

\put(0.1,4.7){\mybox{$0$}} \put(5,4.7){\mybox{$\J \times_{T(M/G)}
T(T(M/G))$}} \put(9.9,4.7){\mybox{$\J \times_{T(M/G)} T(T(M/G))$}}

\end{picture}

By construction $\vlift\Omega\circ\mech\Omega = \mech\Omega$, so
$\Omega^{\J}=\mech\Omega\circ\vlift\Omega$ satisfies
$(\Omega^{\J})^2=\Omega^{\J}$ as it should. We have
\[
\Omega^{\J}(\clift{\tilde E}_A) = \clift{\tilde E}_A,\qquad
\Omega^{\J}(\Hb_\alpha) = 0,\qquad
\Omega^{\J}(\clift{\bar X}_i)=0,\qquad
\Omega^{\J}(\vlift{\bar X}_i)=0.
\]
The tensor field $\Omega^{\J}$ is $G_\mu$-invariant  since both of
the tensors of which it is composed are $G_\mu$-invariant;
$\Omega^{\J}$ therefore defines a principal $G_\mu$-connection.

Note that to define the mechanical connection we do not need a
principal connection $\omega$ on $M\to M/G$ (though we may use one
in calculations). If such a connection is available then we can use
either  $\mech\Omega$ or $\Omega^{\J}$ for the reconstruction
method.

The connection $\Omega^{\J}$ is clearly different from $\mech\Omega$
in general.  We can also decompose $\mech\Omega$ into two connections,
in accordance with the short exact sequences in the diagram.  The
splitting $\vlift{\Omega}_0$ of the middle vertical line, similar to
the connection $\vlift\Omega$ of $\Omega^{\J}$, can be defined by
saying that a vector field $W$ is horizontal if
$g(\tilde\xi,\tau_*W)=0$ for all $\xi\in\g$ (not just for
$\xi\in\g_\mu$).  For this connection we have \[
\vlift{\Omega}_0(\clift{\bar E}_a) = \clift{\bar E}_a,\qquad
\vlift{\Omega}_0( \clift{\hat X}_i)= 0,\qquad \vlift\Omega(
\vlift{\bar X}_i) =0,
\]
where the vector fields $\clift{\hat X}_i$ are exactly those that we
have encountered in Section~2.

To end this section we consider the decomposition of the restriction
of the Euler-Lagrange field $\Gamma$ to $N_\mu$ into its vertical and
horizontal parts with respect to the two connections.

Let us introduce coordinates $(x^i,\theta^a)$ on $M$ such that the
orbits of $G$, or in other words the fibres of $M\to M/G$, are given
by $x^i=\mbox{constant}$; the $x^i$ may therefore be regarded as
coordinates on $M/G$. As before, we will use as fibre coordinates the
quasi-velocities $(v^i,v^a)$ with respect to the standard basis
$\{X_i,{\tilde E}_a\}$.  The non-singularity of $(g_{ab})$ ensures
that, at least locally, we can rewrite the relation $p_a=\mu_a$ for
the injection $\iota:\J\to TM$ in the form $v^a =
\iota^a(x^i,\theta^a,v^i)$, for certain functions $\iota^a$ of the
specified variables.  The restriction of the Euler-Lagrange field to
$\J$ is
\begin{eqnarray*}
\Gamma & = &
\iota^a\clift{\bar{E}_a}+ v^i\clift{\bar{X}_i}+(\Gamma^i\circ
\iota) \vlift{\bar{X}_i}\\
&=&(\iota^A + \Upsilon^A_\alpha  \iota^\alpha)
\clift{\tilde{E}_A}+ {\iota}^\alpha {\Hb}_\alpha + {v}^i
\clift{{\bar X}_i} + (\Gamma^i\circ \iota) \vlift{{\bar X}_i}\\
&=& (\iota^A + \Upsilon^A_\alpha \iota^\alpha + \Upsilon^A_i v^i)
\clift{\tilde{E}_A}+  {\iota}^\alpha {\Hb}_\alpha + {v}^i
\hlift{{\bar X}_i} + (\Gamma^i\circ \iota) \vlift{{\bar X}_i}.
\end{eqnarray*}
The vertical part of $\Gamma$ with respect to the mechanical
connection $\mech\Omega$ is
$(\iota^A + \Upsilon^A_\alpha\iota^\alpha + \Upsilon^A_i v^i)\clift{\tilde{E}_A}$,
and with respect to the vertical lift connection $\Omega^{\J}$ it is
$(\iota^A + \Upsilon^A_\alpha \iota^\alpha)\clift{\tilde{E}_A}$.

Note that neither of the current decompositions of $\Gamma$ coincides
with the one we had towards the end of Section~3, which we should now write
$\Gamma=(\iota^a+B^a_iv^i)\clift{{\bar E}_a}+\Gamma_0$.  The reason is
that this last decomposition is only partial, in the sense that the
vector field $\Gamma_0$ is the horizontal part of $\Gamma$ with
respect to the connection $\vlift{\Omega}_0$; it is the horizontal
lift of a section of the pullback bundle $\J \times_{T(M/G)}
T(T(M/G))$, not a vector field on $T(M/G)$, and this section is only a
part of the data required for the reconstruction method.

\section{The reduced vector field}

A principal connection is all we need to reconstruct integral curves
of an invariant vector field from those of its reduced vector field.
We next examine the latter.

\subsection{The Abelian case}

Before embarking on the more general case, it is instructive to see
what happens if the symmetry group $G$ happens to be Abelian, i.e.\
when $C^c_{ab}=0$.  Then as we pointed out earlier for the case of a
simple mechanical system with Abelian symmetry group, $\g_\mu=\g$
and any level set $p_a=\mu_a$ is invariant under the whole group
$G$.  In fact, under the assumption that $p_a=\mu_a$ can be solved
locally in the form $v^a=\iota^a$, $N_\mu/G$ can be interpreted as
$T(M/G)$, with coordinates $(x^i,v^i)$, where the $x^i$ are
coordinates on $M/G$ and the $v^i$ the corresponding fibre
coordinates (no longer {\em quasi}-velocities).  In this case there
are no `$E_\alpha$'-vectors and $\clift{\bar E}_a = \clift{\tilde
E}_a$ for all $a$.

The restriction of the Euler-Lagrange field to $N_\mu$, given here
by
\begin{eqnarray*}\Gamma &=& \iota^a
\clift{\tilde E}_a + v^i \clift{\bar X}_i+ (\Gamma^i\circ\iota)
\vlift{\bar X}_i \\ &=& (\iota^a - g^{ab}g_{bi}v^i)
\clift{\tilde E}_a + v^i \hlift{\bar X}_i+ (\Gamma^i\circ\iota)
\vlift{\bar X}_i,
\end{eqnarray*}
is now also $G$-invariant.  As a consequence, the coefficients
$\Gamma^i\circ\iota$ do not depend on the group coordinates
$\theta^a$ but only on the coordinates $(x^i,v^i)$ of $T(M/G)$.  In
fact the vector fields $v^i\clift{\bar X}_i+
(\Gamma^i\circ\iota)\vlift{\bar X}_i$ (the $\Omega^{\J}$-horizontal
part of $\Gamma$) and $v^i \hlift{\bar X}_i+
(\Gamma^i\circ\iota)\vlift{\bar X}_i$ (the $\mech\Omega$-horizontal
part of $\Gamma$) both reduce to the same vector field on $N_\mu/G$,
which in this case is exactly a second-order differential equation
field on $T(M/G)$,
\[
\check\Gamma = v^i\fpd{}{x^i} + \Gamma^i(x,v)\fpd{}{v^i}.
\]
The integral curves of this reduced vector field are the solutions
of the equations ${\ddot x}^i=\Gamma^i(x,\dot x)$ (with
$v^i=\dot{x}^i$) and, from the introduction, we know that these are
equivalent to the equations
\[
\frac{d}{dt}\left(\fpd{\R}{v^i}\right)-\fpd{\R}{x^i}= -
B^a_{ij}\pi_a{\dot x}^j.
\]

\subsection{The non-Abelian case}

In the general case of a non-Abelian symmetry group we should not
expect that the equations for $x^i$ will be completely decoupled from
all coordinates $\theta^a$.  Indeed, in that case the vector field
$\Gamma$ reduces to a vector field on $N_\mu/G_\mu$.  This manifold can
locally be identified with $M/G_\mu\times T(M/G)$, so the equations
for the integral curves of the reduced vector field will depend also
on the coordinates of $M/G_\mu$.

To give a local expression of the reduced vector field
$\check\Gamma$ we need to introduce a basis for
$\vectorfields{\J/G_\mu}$.  The bracket relations $[\clift{{\tilde
E}_A}, \clift{{\bar X}_i}]=0$ and $[\clift{{\tilde E}_A},
\vlift{{\bar X}_i}]=0$ show that $\clift{\bar X}_i$ and $\vlift{\bar
X}_i$ are $G_\mu$-invariant vector fields on $\J$; they project
therefore onto vector fields $\clift{\check X}_i$ and $\vlift{\check
X}_i$ on $\J/G_\mu$.  The invariance of the Hessian $g$ amounts for
its coefficients to $\clift{\tilde{E}_a}(g_{ij})=0$,
$\clift{\tilde{E}_a}(g_{bc})+C^d_{ab}g_{cd}+C^d_{ac}g_{bd}=0$, and $
\clift{\tilde{E}_a}(g_{ib})+C^c_{ab}g_{ic}=0$.  From this
\begin{eqnarray*}
\clift{\tilde E}_A(\Upsilon^B_i) &=& C^B_{AC}\Upsilon^C_i,\\
\clift{\tilde E}_A(\Upsilon^B_\alpha) &=&
C^B_{AC}\Upsilon^C_\alpha
-C^\beta_{A\alpha}\Upsilon^B_\beta-C^B_{A\alpha}
\end{eqnarray*}
(where we have taken into account the fact that in the current basis
$C^\gamma_{AB}=0$).  It is now easy to see that the vector fields
$\hlift{\bar X}_i =\clift{\bar X}_i - \Upsilon_i^A \clift{\tilde E}_A$
are also $G_\mu$-invariant.  In fact, since they differ from
$\clift{\bar X}_i$ only in a part that is vertical with respect to the
bundle projection $N_\mu\to N_\mu/G_\mu$, they project onto the same
vector fields $\clift{\check X}_i$ on $N_\mu/G_\mu$.

The vector fields $\Hb_\alpha$ are not invariant:\ in fact
$[\clift{{\tilde E}_A}, \Hb_\alpha] = -C_{A\alpha}^\beta \Hb_\beta$.
To obtain a complete basis for $\vectorfields{N_\mu/G_\mu}$, we need
to replace the vector fields $\{\Hb_\alpha\}$ by $G_\mu$-invariant
vector fields. To do so, we will consider the $G$-invariant vector
fields ${\hat E}_a={\mathcal A}_a^b {\tilde E}_b$ on $M$ that we
introduced in Section~2. Let $({\mathcal A}^\alpha_\beta)$ be the
coefficients we find in the relation $ {\hat E}_\alpha = {\mathcal
A}_\alpha^\beta{\tilde E}_\beta + {\mathcal A}_\alpha^B {\tilde
E}_B$. The vector fields
\[
{\Hh}_\alpha = {\mathcal A}_\alpha^\beta \Hb_\beta
\]
are tangent to the level set $N_\mu$ and horizontal. Given that
$C^\beta_{AB}=0$, it easily follows from the relation
${\tilde E}_A({\mathcal A}_\alpha^\beta)=
C^\beta_{A\gamma}{\mathcal A}^\gamma_\alpha$ that these vector fields are
$G_\mu$-invariant:
\[
[\clift{{\tilde E}_A}, {\Hh}_\alpha]=
{\tilde E}_A({\mathcal A}_\alpha^\beta) \Hb_\beta
- {\mathcal A}_\alpha^\beta C_{A\beta}^\gamma \Hb_\gamma=0.
\]
They project therefore onto vector fields ${\F}_\alpha$ on
$\J/G_\mu$. To conclude, the set $\{ \clift{\check
X}_i,\vlift{\check X}_i, {\F}_\alpha \}$ defines the basis for
$\vectorfields{\J/G_\mu}$ we were looking for.

We denote by $({\bar {\mathcal A}}^\alpha_\beta)$ the matrix inverse
to $ ({\mathcal A}^\alpha_\beta)$ and set $\Phi^A=\iota^A +
\Upsilon^A_\alpha  \iota^\alpha$, $\submech{\Phi^A}=\iota^A +
\Upsilon^A_\alpha \iota^\alpha + \Upsilon^A_i v^i$ and $\Psi^\alpha
= {\bar {\mathcal A}}^\alpha_\beta {\iota}^\beta$. Then $\Gamma$
takes the form
\begin{eqnarray*}
\Gamma & = &  \Phi^A \clift{\tilde{E}_A}+ \Psi^\alpha {\Hh}_\alpha +
{v}^i \clift{{\bar X}_i} + (\Gamma^i\circ \iota) \vlift{{\bar
X}_i},\\ & = &  \submech{\Phi^A}\clift{\tilde{E}_A}+ \Psi^\alpha
{\Hh}_\alpha + {v}^i \hlift{{\bar X}_i} + (\Gamma^i\circ \iota)
\vlift{{\bar X}_i},
\end{eqnarray*}
where, as before, the first term is the vertical part of $\Gamma$
with respect to the vertical lift connection $\Omega^{\J}$ (in the
first place) and the mechanical connection $\mech\Omega$ (in the
second). Obviously $v^i$ and $(\Gamma^i \circ \iota)$ are
$G_\mu$-invariant functions.  To see that $\Psi^\alpha$ is also
$G_\mu$-invariant, recall that $C^\alpha_{AB}=0$ and $\clift{\tilde
E}_A (\iota^\beta) = C^\beta_{A\gamma}\iota^\gamma$, and observe
that $\clift{\tilde E}_A ({\bar {\mathcal A}}^\beta_\delta)
{\mathcal A}^\gamma_\beta = -{\bar {\mathcal A}}^\beta_\delta
\clift{\tilde E}_A ( {\mathcal A}^\gamma_\beta)= -{\bar {\mathcal
A}}^\beta_\delta {\mathcal A}^\alpha_\beta C^\gamma_{A\alpha}=
-C^\gamma_{A\delta}$. Therefore
\[
\clift{\tilde E}_A (\Psi^\alpha)=\clift{\tilde E}_A ( {\bar
{\mathcal A}}^\alpha_\beta) {\iota}^\beta +{\bar {\mathcal
A}}^\alpha_\beta \clift{\tilde E}_A ( {\iota}^\beta) =  -{\bar
{\mathcal A}}^\alpha_\delta C^\delta_{A\beta}\iota^\beta + {\bar
{\mathcal A}}^\alpha_\beta C^\beta_{A\gamma}\iota^\gamma=0.
\]
We conclude that $v^i$, $(\Gamma^i\circ \iota)$ and $\Psi^\alpha$
can all be regarded as functions on $\J/G_\mu$. The horizontal part
of $\Gamma$, for both connections,  can thus be interpreted as the
horizontal lift of the reduced vector field
\[
\check\Gamma = \Psi^\alpha {\F}_\alpha + {v}^i \clift{{\check X}_i}
+ (\Gamma^i\circ \iota) \vlift{{\check X}_i}
\]
 on $\J/G_\mu$.

For completeness, we point out that it follows from the relations
$\clift{\tilde E}_A(\Upsilon^B_\alpha)=C^B_{AC}\Upsilon^C_\alpha
-C^\beta_{A\alpha}\Upsilon^B_\beta-C^B_{A\alpha}$ and $\clift{\tilde
E}_A(\iota^B) = C^B_{AC}\iota^C+C^B_{A\gamma}\iota^\gamma$ that the
coefficients $\Phi^A$ and $\submech{\Phi^A}$ satisfy
\[
\clift{\tilde E}_A(\Phi^B) = C^B_{AC}\Phi^C, \qquad
\clift{\tilde E}_A(\submech{\Phi^B}) = C^B_{AC}\submech{\Phi^C}.
\]
This shows that they can be interpreted as the coefficients of
$\g_\mu$-valued functions $\Phi$ and $\submech\Phi$ on $\J$ satisfying
$\Phi\circ{\psi_g}^{N_\mu} = \ad_{g}\Phi$ for $g\in G_\mu$ (and
similarly for $\submech\Phi$), where $\psi^{N_\mu}$ denotes the
$G_\mu$-action on $\J$ (see \cite{Paper1}).

We now give a coordinate expression for the reduced vector field.
From here on we will use coordinates
$(\theta^a)=(\theta^A,\theta^\alpha)$ such that the fibres of $G\to
G/G_\mu$ are given by $\theta^\alpha=\mbox{constant}$.  With this
assumption, there are functions $K^a_b$ on $M$ such that
\[
{\tilde E}_A = K^B_A \fpd{}{\theta^B}, \qquad {\tilde E}_\alpha =
K^B_\alpha \fpd{}{\theta^B} + K^\beta_\alpha \fpd{}{\theta^\beta}.
\]
We also introduce the functions $\Lambda^b_i$ for which
\[
X_i =\fpd{}{x^i} - {\Lambda}^b_i \fpd{}{\theta^b},
\]
as before.

By interpreting $\J/G_\mu$ locally as $ M/G_\mu\times_{M/G} T(M/G)$,
we see that a point of $\J/G_\mu$ has coordinates
$(x^i,\theta^\alpha,v^i)$.  Because of their $G_\mu$-invariance, the
functions $\Gamma^i\circ \iota$ and $\Psi^\alpha$ are independent of
the variables $\theta^A$.  Let $\pi^{N_\mu}$ be the projection $N_\mu
\to N_\mu/G_\mu$; then for any invariant function $F$ on $N_\mu$ there
is a function $f$ on $N_\mu/G_\mu$ such that $F=f\circ\pi^{N_\mu}$.
Then for all invariant vector fields $X$ on $N_\mu$, and their
reductions $\check X$ to vector fields on $N_\mu/G_\mu$, we have
\[
X(F) = \check X(f)\circ\pi^{N_\mu}.
\]
We will apply this property to the vector fields $\hlift{\hat
E}_\alpha$, $\clift{\bar X}_i$ and $\vlift{\bar X}_i$ and the
invariant functions $x^i$, $v^i$ and $\theta^\alpha$. Keeping in
mind that for any vector field $Z$, function $f$  and 1-form
$\theta$ on $M$,
\[
\clift Z (\tau^*f) = \tau^*Z(f), \quad \vlift Z (\tau^*f)=0,\quad
\clift{Z}(\vec{\theta})=\overrightarrow{\lie{Z}\theta}, \quad
\vlift{Z}(\vec{\theta})=\tau^*\theta(Z).
\]
where $\vec\theta$ stands for the fibre-linear function on $TM$
defined by the 1-form $\theta$, and $\tau$ is the tangent projection
$TM \to M$, we find that
\[
\begin{array}{lclcl}
\clift{\bar X}_i (x^j) = \delta_i^j,& & \clift{\bar X}_i
(\theta^\beta)=-\Lambda_i^\beta, & & \clift{\bar X}_i (v^j) = 0,\\
 \vlift{\bar X}_i (x^j) =0,&  & \vlift{\bar X}_i
(\theta^\beta)=0, & & \vlift{\bar X}_i (v^j) = \delta_i^j, \\
\hlift{\hat E}_\alpha (x^j) = 0, & & \hlift{\hat E}_\alpha
(\theta^\beta) = {\mathcal A}^\gamma_\alpha K^\beta_\gamma, & &
\hlift{\hat E}_\alpha (v^j) = 0,
\end{array}
\]
from which it follows immediately that
\[
\clift{{\check X}_i}= \fpd{}{x^i} - {\Lambda}^\alpha_i
\fpd{}{\theta^\alpha}, \qquad \vlift{{\check X}_i} =
\fpd{}{v^i},\qquad {\F}_\alpha = {{\mathcal A}}^\gamma_\alpha
K^\beta_\gamma\fpd{}{\theta^\beta}
\]
and
\[
\check\Gamma= (\iota^\beta K^\alpha_\beta - v^i
{\Lambda}^\alpha_i)\fpd{}{\theta^\alpha} + {v}^i \fpd{}{x^i} +
(\Gamma^i\circ \iota) \fpd{}{v^i}.
\]
The equations that determine the integral curves
${\check v}(t)=(x^i(t),\theta^\alpha(t),v^i(t))$ of the reduced vector
field $\check\Gamma$ are therefore the coupled set
\[
\left\{\begin{array}{lll}
\ddot{x}^i &=& \Gamma^i\circ\iota, \\
\dot{\theta}^\alpha &=& \iota^\beta K^\alpha_\beta -
{v}^i{\Lambda}^\alpha_i.\end{array} \right.
\]
One can easily convince oneself that the right hand side of the
equation for ${\dot \theta}^\alpha$ is indeed independent of the
variables $\theta^A$:\ by considering the coefficients of
$\partial/\partial\theta^\alpha$ in $[{\tilde E}_A,{\tilde
E}_\beta]=-C_{A\beta}^c{\tilde E}_c$ we find that ${\tilde
E}_A(K^\alpha_\beta)=-C^\gamma_{A\beta}K^\alpha_\gamma$; since also
$\clift{\tilde E}_A(\iota^\beta)=C^\beta_{A\gamma}\iota^\gamma$ and
${\tilde E}_A({\Lambda}^\beta_i)=0$, it follows easily that
$\clift{\tilde E}_A(\iota^\beta K^\alpha_\beta -
v^i{\Lambda}^\alpha_i)=0$, as claimed.

The functions $\Gamma^i\circ\iota$ on the right-hand side of the
equation for the $\ddot{x}^i$ can be determined from the generalized
Routh equations of Section~3,
\[
\Gamma(\vlift{\bar{X}_i}(\R^\mu))-\clift{\bar{X}_i}(\R^\mu)=-\mu_aR^a_{ij}v^j.
\]
Since $\R^\mu=\R\circ\iota$ is $G_\mu$-invariant, so also are
$\vlift{\bar{X}_i}(\R^\mu)$ and $\clift{\bar{X}_i}(\R^\mu)$.  Recall
that the $R^a_{ij}$ are functions on $M$, determined by
$[X_i,X_j]=R^a_{ij}\tilde{E}_a$. Thus since $[\tilde{E}_a,X_i]=0$,
\[
(\tilde{E}_a(R^b_{ij})-R^c_{ij}C_{ac}^b)\tilde{E}_b=0.
\]
But $C_{Ac}^b\mu_b=0$, and so
\[
\clift{\tilde{E}_A}(R^b_{ij}\mu_b)=
\tilde{E}_A(R^b_{ij})\mu_b=
R^c_{ij}C_{Ac}^b\mu_b=0.
\]
It follows that the term $\mu_aR^a_{ij}v^j$ is $G_\mu$-invariant.  The
generalized Routh equations therefore pass to the quotient
$N_\mu/G_\mu$, and take the reduced form
\[
{\check\Gamma}(\vlift{\check{X}_i}(\R^\mu))-\clift{\check{X}_i}(\R^\mu)=
-\mu_aR^a_{ij}v^j.
\]
Following \cite{MRS}, we will call these reduced equations the
Lagrange-Routh equations.  Under the regularity assumptions we have
adopted throughout, the function-valued matrix $(\vlift{\bar X}_i
\vlift{\bar X}_j(\R^\mu))$ is non-singular and the coefficients
$\Gamma^i\circ\iota$, now interpreted as functions on $\J/G_\mu$,
can be determined from the Lagrange-Routh equations.  In the current
coordinate system the equations become
\[
\frac{d}{dt}\left(\fpd{\R^\mu}{v^i}\right)-\fpd{\R^\mu}{x^i}
=-\mu_aR^a_{ij}v^j - \Lambda^\alpha_i \fpd{\R^\mu}{\theta^\alpha}.
\]
Given a reduced solution $\check
v(t)=(x^i(t),\theta^\alpha(t),v^i(t))\in N_\mu/G_\mu$, we can apply
the method of reconstruction using either one of the connections
$\mech\Omega$ and $\Omega^{N_\mu}$ to recover a complete solution
$v(t)=(x^i(t),\theta^A(t),\theta^\alpha(t),v^i(t))\in N_\mu$ of the
Lagrangian system.  The examples discussed in Section~\ref{ex} will
make it clear how this method works in practice.

\section{Simple mechanical systems}

In this section we reconcile our results with those for
the case of a simple mechanical system to be found elsewhere in the
literature.

A simple mechanical system is one whose Lagrangian is of the form
$L=T-V$ where $T$ is a kinetic energy function, defined by a
Riemannian metric $g$ on $M$, and $V$ is a function on $M$, the
potential energy.  The symmetry group $G$ consists of those
isometries of $g$ which leave $V$ invariant.  We define a connection
on $M\to M/G$ by taking for horizontal subspaces the orthogonal
complements to the tangent spaces to the fibres; it is this
connection that is usually called the mechanical connection (for a
simple mechanical system).  We write
$g_{ab}=g(\tilde{E}_a,\tilde{E}_b)$, $g_{ij}=g(X_i,X_j)$; by
assumption, $g(\tilde{E}_a,X_i)=0$.  Then in terms of
quasi-velocities,
\[
L(m,v)=\onehalf\left(g_{ij}(m)v^iv^j+g_{ab}(m)v^av^b\right)-V(m),
\]
and $\vlift{\tilde{E}_a}\vlift{\tilde{E}_b}(L)=g_{ab}$ etc., so the
notation is consistent with what has gone before.  Note that since
we assume that $g$ is Riemannian and therefore positive-definite it
is automatic that $L$ is regular and that the matrices $(g_{ab}(m))$
and $(g_{ij}(m))$ are both non-singular for all $m$; in particular
we don't need to make the separate assumption that $(g_{ab})$ is
non-singular. Considered as defining a map $M\to \g^*\odot\g^*$,
$(g_{ab})$ is called the locked inertia tensor.  The isometry
condition gives
\[
\tilde{E}_a(g_{bc})+C_{ab}^dg_{cd}+C_{ac}^dg_{bd}=0,\quad
\tilde{E}_a(g_{ij})=0.
\]
The first of these is the differential version of the equivariance
property of the locked inertia tensor with respect to the action of
$G$ on $M$ and the coadjoint action of $G$ on $\g^*\odot\g^*$.  The
second shows that $g_{ij}$ may be considered as a function on $M/G$.

The momentum is given simply by $p_a(m,v)=g_{ab}(m)v^b$.  On any
level set $N_\mu$, where $p_a=\mu_a$, we can solve explicitly for
the $v^a$ to obtain $v^a=g^{ab}\mu_b$.

The Routhian is given by
\[
\R=L-p_av^a=\onehalf g_{ij}v^iv^j-\onehalf g_{ab}v^av^b-V,
\]
and on restriction to $N_\mu$ we obtain
\[
\R^\mu=\onehalf g_{ij}v^iv^j-\left(V+\onehalf
g^{ab}\mu_a\mu_b\right).
\]
The quantity $V+\onehalf g^{ab}\mu_a\mu_b$ is the so-called amended
potential \cite{Simo} and the term $C_\mu=\onehalf g^{ab}\mu_a\mu_b$
is called the `amendment' in \cite{MRS}.  Both functions on $M$ are
$G_\mu$-invariant: one easily verifies that
$\tilde{E}_a(C_\mu)=g^{bc}C^d_{ab}\mu_c\mu_d$, so in particular for
$a=A$ we get $\tilde{E}_A(C_\mu)=0$.

Note that by the choice of connection $B^a_i=0$; we have $\clift{\hat
X}_i=\clift{\bar X}_i$, and the generalized Routh equations are
\[
\Gamma_0(\vlift{\bar{X}_i}(\R^\mu))
-\clift{\bar{X}_i}(\R^\mu)=-\mu_aR^a_{ij}v^j.
\]
This equation is the analog in our framework of the one in
Corollary~III.8 of \cite{MRS}.  We have shown in the previous section
that it reduces to the Lagrange-Routh equations
\[
{\check\Gamma}(\vlift{\check{X}_i}(\R^\mu))-\clift{\check{X}_i}(\R^\mu)=
-\mu_aR^a_{ij}v^j,
\]
which for consistency should be supplemented by the equation that
determines the variables $\theta^\alpha$.  As we pointed out earlier,
the latter is actually just the expression for genuine velocity
components $\dot{\theta}^\alpha$ in terms of quasi-velocities,
supplemented by the constraint $v^\alpha=\iota^\alpha$ which for a
simple mechanical system takes the form
$\iota^\alpha= g^{\alpha a}\mu_a$.

We can split the reduced Routhian $\R^\mu$ in a Lagrangian part
$\goth{L}=\onehalf g_{ij}v^iv^j-V$ and the reduced amendment
$\goth{C}_\mu$. Since the quasi-velocities $v^a$ do not appear in
the expression of $\goth{L}$, it can formally be interpreted as a
function on $T(Q/G)$. The reduced amendment is a function on
$Q/G_\mu$. We can now rewrite the Lagrange-Routh equations as
\[
{\check\Gamma}(\vlift{\check{X}_i}(\goth{L}))-\clift{\check{X}_i}(\goth{L})=
-\mu_aR^a_{ij}v^j + \check{X}_i(\goth{C}_\mu);
\]
in coordinates
\[
\frac{d}{dt}\left(\fpd{\goth{L}}{v^i}\right)-\fpd{\goth{L}}{x^i}=
-\mu_aR^a_{ij}v^j +
\left(\vf{x^i}-\Lambda^\alpha_i\vf{\theta_\alpha}\right)(\goth{C}_\mu).
\]
This equation is only one out of two equations that appear in
Theorem~III.14 in \cite{MRS}, the theorem that states the reduced
equations obtained by following a variational approach to Routh's
procedure.  We leave it to the reader to verify that the second
equation, in its form (III.37), is in fact
\[
v^\alpha C^a_{\beta\alpha}\mu_a
= g^{\alpha b}\mu_b C^a_{\beta\alpha}\mu_a.
\]
Since $v^\alpha=g^{\alpha b}\mu_b$, this is obviously an identity
from the current point of view; it certainly cannot be used to
determine ${\dot\theta}^\alpha$ in terms of the other variables, and
without this information the equations are incomplete.  In this
respect, therefore, our reduction results are an improvement on
those in \cite{MRS}.

Let us now check, in the case where the configuration space $M$ is of
the form $S\times G$, for an Abelian symmetry group ($C^c_{ab}=0$),
and a Lagrangian of the form
\[
L(x,\theta,\dot x,\dot\theta)= \onehalf k_{ij}(x){\dot x}^i{\dot
x}^j + k_{ia}(x) {\dot x}^i{\dot \theta}^a + \onehalf
k_{ab}(x){\dot \theta}^a{\dot \theta}^b-V(x),
\]
that the reduced equations above coincide with those in the
introduction. We set
\[
{\tilde E}_a = K_a^b\fpd{}{\theta^b}
\]
where the $K^a_b$ are independent of the $\theta^a$ since we are
dealing with the Abelian case.  In general, horizontal vector fields
take the form
\[
X_i= \fpd{}{x^i} - \Lambda^a_i \fpd{}{\theta^a}.
\]
The quasi-velocities adapted to the connection are therefore given, as
before, by $v^i={\dot x}^i$ and $K^a_bv^b={\dot\theta}^a + \Lambda^a_i
{\dot x}^i$.

Given that in this case $\tilde{E}_a(g_{bc})=0$ and
$\tilde{E}_a(g_{ij})=0$, all coefficients of the metric can be
interpreted as functions on $M/G=S$, and they depend only on the
variables $x^i$.  The use of the mechanical connection entails that
$g_{ai}=0$.  When expressed in terms of the coordinates $({\dot
x}^i,{\dot\theta}^a)$, this property fixes the connection coefficients
to be of the form $\Lambda^a_i=k^{ab}k_{ib}$ and the remaining
coefficients of the metric to be $g_{ab}=k_{cd}K^c_aK^d_b$ and
$g_{ij}=k_{ij}-k^{ab}k_{ia}k_{jb}$.  The expression for $\R^\mu$ given
in the introduction now easily follows.  Since $\R^\mu$ is a function
only of $x^i$ and $v^i=\dot{x}^i$ we get
\[
\vlift{\bar{X}_i}(\R^\mu) = \fpd{\R^\mu}{{\dot x}^i} \quad
\mbox{and}\quad \clift{\bar{X}_i}(\R^\mu) = \fpd{\R^\mu}{x^i}.
\]
Moreover,
\[
[X_i,X_j] = B_{ij}^a \fpd{}{\theta^a} =
R_{ij}^a{\tilde E}_a,
\]
and likewise $\mu_a = \vlift{{\tilde E}_a}(L)=K_a^b\pi_b$. So
$\mu_aR^a_{ij} = \pi_a B^a_{ij}$ and the equation from the
introduction follows.

We now return to the general case (of a simple mechanical system)
and consider the reconstruction process.

We continue to use the mechanical connection on $M$.  Since, in the
basis that is adapted to this connection, $g_{ia}=0$, and therefore
$\hlift{\bar X}_i = \clift{\bar X}_i$, the two connections
$\mech\Omega$ and $\Omega^{N_\mu}$ coincide.  We denote the common
connection on $N_\mu$ by $\Omega$.

Let $\check{v}(t)$ be a curve in $N_\mu/G_\mu$ which is an integral
curve of $\check{\Gamma}$, and $\hlift{\check{v}}$ a horizontal lift of
$\check{v}$ to $N_\mu$ (horizontal with respect to $\Omega$). The
reconstruction equation is
\[
\clift{\widetilde{\vartheta(\dot{g})}}=\Omega(\Gamma\circ\hlift{\check{v}})
\]
(where $\vartheta$ here is the Maurer-Cartan form of $G_\mu$); this is (at
each point on the curve $\check{v}$) an equation between vertical
vectors on $N^\mu$, but can and should be thought of as an equation on
$\g_\mu$.  It determines a curve $g(t)$ in $G_\mu$ such that
\[
t\mapsto \psi^{N_\mu}_{g(t)}\hlift{\check{v}}(t)
\]
is an integral curve of $\Gamma$ in $N_\mu$; again,
$\psi^{N_\mu}$ is the action of $G_\mu$ on $N_\mu$.  So far, this
works for an arbitrary Lagrangian.

Now as we showed earlier in general, the vertical part of $\Gamma$
with respect to the vertical lift connection $\Omega^{\J}$ is
$(\iota^A + \Upsilon^A_\alpha \iota^\alpha)\clift{\tilde{E}_A}$.
Thus in the case at hand
\begin{eqnarray*}
\Omega(\Gamma)&=&(\iota^A + \Upsilon^A_\alpha
\iota^\alpha)\clift{\tilde{E}_A}\\
&=&
(\iota^A + G^{AB}g_{B\alpha}\iota^\alpha)\clift{\tilde{E}_A}\\
&=&(g^{Aa}\mu_a + G^{AB}g_{B\alpha}g^{\alpha a}\mu_a)\clift{\tilde{E}_A}\\
&=&(g^{Aa}\mu_a +
G^{AB}(\delta^a_B - g_{BC}g^{Ca})\mu_a)\clift{\tilde{E}_A}
\\&=&(g^{Aa}\mu_a +
G^{AB}\mu_B - g^{Aa}\mu_a)\clift{\tilde{E}_A}\\
&=& G^{AB}\mu_B\clift{\tilde{E}_A}.
\end{eqnarray*}
The first point to note is that the coefficient $G^{AB}\mu_B$
appearing on the right-hand side of the final equation above is a
function on $M$, so that in the right-hand side of the reconstruction
equation the argument $\hlift{\check{v}}$ can be replaced by its
projection into $M$, which is $\tau\circ\iota\circ\hlift{\check{v}}$.

Next, we interpret $G^{AB}\mu_B$ in terms of the locked inertia tensor.
Recall that the locked inertia tensor at $m\in M$ has components $g_{ab}(m)$.
As is the usual practice we consider the locked inertia tensor as a
non-singular symmetric linear map $I(m): \g\to\g^*$.
Now let $j$ be the injection $\g_\mu\to\g$:\ then $\mu_B$ are the
components of $j^*\mu\in\g_\mu^*$, and $g_{AB}(m)$ are the components
of the map $I_\mu(m)=j^*\circ I(m)\circ j$. Then
\[
G^{AB}(m)\mu_BE_A=I^{-1}_\mu(m)(j^*\mu),
\]
a point of $\g_\mu$. So finally the reconstruction equation may be
written
\[
\vartheta(\dot{g}(t))=I^{-1}_\mu(c(t))(j^*\mu),\qquad
c=\tau\circ\iota\circ\hlift{\check{v}}.
\]
This is an equation between curves in $\g_\mu$.

We will now show that this reconstruction equation above is a
particular and simple case of one of the reconstruction equations
appearing in \cite{MRS}.

To do so we must introduce yet another connection, used in \cite{MRS}
and called there the mechanical connection for the $G^\mu$-action.
This is a connection on the principal fibre bundle $M \to M/G_\mu$,
i.e.\ a $G_\mu$-invariant splitting of the short exact sequence
\[
0\to M\times\g_\mu \to TM \to M \times_{M/G_\mu} T(M/G_\mu)\to 0,
\]
If, as before, $\{E_A,E_\alpha\}$ is a basis of $\g$ for which
$\{E_A\}$ is a basis for $\g_\mu$, then the vector fields $X_i$
together with the vector fields ${\tilde E}_\alpha -
G^{AB}g_{A\alpha} {\tilde E}_B$ form a basis for the set of vector
fields which are horizontal with respect to the mechanical
connection for the $G^\mu$-action.  We denote the latter by
$\omega^\mu$. Now $\omega^\mu$ and $\Omega$ are related somewhat as
a connection and its vertical lift:\ in fact (for their projections
onto $\g_\mu$)
\[
\Omega(Z_v) = \omega^\mu(T(\tau\circ\iota) Z_v ), \qquad Z_v \in
T\J.
\]
We note in passing that since $T(\tau\circ\iota)\Gamma(v)=v$ for any
$v\in N_\mu$, we can write the reconstruction equation as
\[
\widetilde{\vartheta(\dot{g}(t))}=\omega^\mu(\hlift{\check{v}}).
\]

The reconstruction equation in \cite{MRS} that we are aiming for is
the third of the four, equation (IV.6).  It seems the
one most relevant to our approach because, as Marsden et al.\ say, in it
they `take the dynamics into account', and this has been our purpose
throughout.  Now equation (IV.6) of \cite{MRS} differs from our
reconstruction equation (expressed in terms of $I_\mu$) by having an
additional term on the right-hand side involving the mechanical
connection for the $G^\mu$-action $\omega^\mu$.  This arises because
the authors start with a more general class of curves on $M$ than we
do.

In order to show that our equation agrees with theirs we first
show that the curve $c=\tau\circ\iota\circ\hlift{\check{v}}$ in $M$ is
$\omega^\mu$-horizontal; the extra term in their equation is
therefore zero in our case. By evaluating $\omega^\mu$ on the tangent
to $c$ and using the relation between $\omega^\mu$ and $\Omega$ we have
\[
\omega^\mu(\dot{c})=
\omega^\mu(T(\tau\circ\iota)\dot{\hlift{\check{v}}})
=\Omega(\dot{\hlift{\check{v}}})=0
\]
because $\hlift{\check{v}}$ is $\Omega$-horizontal.  So our
reconstruction equation formally agrees with equation (IV.6) of
Marsden et al., when we take the starting curve on $M$ to be $c$:\ it
is the particular case of that equation in which the curve on $M$ is
horizontal with respect to the $G_\mu$ mechanical connection.

To finish the story we must also take into account the fact that
equation (IV.6) of \cite{MRS} is presented as an equation for the
reconstruction of a base integral curve of $\Gamma$, with momentum
$\mu$, from another suitable curve on $M$, whereas our reconstruction
equation gives an integral curve of $\Gamma$ on $N^\mu$.  But there is
no real discrepancy here, because $\Gamma$ is a second-order
differential equation field and so knowing its base integral curves is
equivalent to knowing its integral curves.  Let us spell this out in
detail.  We know that if $t\mapsto g(t)$ is a solution of our
reconstruction equation then
\[
t\mapsto \psi^{N_\mu}_{g(t)}\hlift{\check{v}}(t)
\]
is an integral curve of $\Gamma$ in $N_\mu$. The corresponding base
integral curve is
\[
t\mapsto \tau(\iota(\psi^{N_\mu}_{g(t)}\hlift{\check{v}}(t))).
\]
But
\[
\tau\circ\iota\circ\psi^{N_\mu}_g=
\tau\circ\psi^{TM}_g\circ\iota=
\psi^M_g\circ\tau\circ\iota,
\]
so the curve $t\mapsto \psi^M_{g(t)}c(t)$ is a base integral curve of
$\Gamma$.  Thus the same curve in $G_\mu$ determines an integral curve
of $\Gamma$ in $N_\mu$ (by its action on $\hlift{\check{v}}$) and the
corresponding base integral curve (by its action on
$c=\tau\circ\iota\circ\hlift{\check{v}}$, the projection of
$\hlift{\check{v}}$ to $M$).

\section{Illustrative examples}\label{ex}

We give two examples. In the first we derive Wong's equations using
our methods. This example is intended to illustrate the Routhian
approach in a case of some physical interest; however, we do not
pursue the calculations as far as the consideration of the isotropy
algebra and reconstruction. These matters are illustrated in the
second example, which is more specific and more detailed, if
somewhat more artificial.

\subsection{Wong's equations}

We discuss the generalized Routh equations for the geodesic field of a
Riemannian manifold on which a group $G$ acts freely and properly to
the left as isometries, and where the vertical part of the metric
(that is, its restriction to the fibres of $\pM:M\to M/G$) comes from
a bi-invariant metric on $G$.  The reduced equations in such a case
are known as Wong's equations~\cite{Cendra,Mont}.

This is of course an example of a simple mechanical system, with
$V=0$; we therefore adopt the notation of Section~7, and we will use
the mechanical connection.  In order to utilise conveniently the
assumption about the vertical part of the metric $g$, we will need
symbols for the components of $g$ with respect to the invariant vector
fields $\hat{E}_a$ introduced in Section~2; we write
\[
h_{ab}=g(\hat{E}_a,\hat{E}_b)=\A_a^c\A_b^dg_{cd}.
\]
Since both $h_{ab}$ and $g_{ij}$ are $G$-invariant functions, they
pass to the quotient.  In particular, the $g_{ij}$ are the components
with respect to the coordinate fields of a metric on $M/G$, the
reduced metric; we denote by $\conn kij$ its Christoffel symbols.

The further assumption about the vertical part of the metric has the
following implications.  It means in the first place that
$\lie{\hat{E}_c}g(\hat{E}_a,\hat{E}_b)=0$ (as well as
$\lie{\tilde{E}_c}g(\hat{E}_a,\hat{E}_b)=0$).  Taking into account the
bracket relations $[\hat{E}_a,\hat{E}_b]=C^c_{ab}\hat{E}_c$, we find
that the $h_{ab}$ must satisfy $h_{ad}C^d_{bc}+h_{bd}C^d_{ac}=0$.  It
is implicit in our choice of an invariant basis that we are working in
a local trivialization of $M\to M/G$.  Then the $h_{ab}$ are functions
on the $G$ factor, so must be independent of the coordinates $x^i$ on
$M/G$, which is to say that they must be constants.  Moreover,
$\tilde{E}_a$, $\hat{E}_a$ and $\A^b_a$ are all objects defined on the
$G$ factor, so are independent of the $x^i$.  We may write
\[
X_i=\vf{x^i}-\gamma_i^a\hat{E}_a
\]
for some coefficients $\gamma_i^a$ which are clearly $G$-invariant;
moreover $[X_i,\hat{E}_a]=\gamma_i^cC_{ac}^b\hat{E}_b$.  We set
$\gamma_i^cC_{ac}^b=\gamma^b_{ia}$; then
$h_{ac}\gamma^c_{ib}+h_{bc}\gamma^c_{ib}=0$.

We are interested in the geodesic field of the Riemannian metric $g$.
The geodesic equations may be derived from the Lagrangian
\[
L=\onehalf g_{\alpha\beta}u^\alpha u^\beta
=\onehalf g_{ij}v^iv^j+\onehalf g_{ab}v^av^b
=\onehalf g_{ij}v^iv^j+\onehalf h_{ab}w^aw^b,
\]
where the $w^a$ are quasi-velocities relative to the $\hat{E}_a$; we
have $\A^a_bw^b=v^a$.
The momentum is given by $p_a=g_{ab}v^b=\bar{\A}_a^ch_{bc}w^c$, where
$(\bar{\A}_a^b)$ is the matrix inverse to $(\A_a^b)$.
The Routhian is
\[
\R=\onehalf g_{ij}v^iv^j-\onehalf g^{ab}p_ap_b.
\]
It is easy to see that
$\vlift{\bar{X}_i}(\R)=g_{ij}v^j$. The calculation of
$\clift{\bar{X}_i}(\R)$ reduces to the calculation of $X_i(g_{ij})$
and $X_i(g^{ab})$. The first is straightforward. For the second, we
note that $g_{ab}=\bar{\A}_a^c\bar{\A}_b^dh_{cd}$; since the
right-hand side is independent of the $x^i$, so is $g_{ab}$, and so
equally is $g^{ab}$.  It follows that
\[
\clift{\bar{X}_i}(\R)=
\onehalf\fpd{g_{jk}}{x^i}v^jv^k-\onehalf\gamma^c_i\hat{E}_c(g^{ab})p_ap_b.
\]
Now $\hat{E}_c(g^{ab})=-\A_c^d(g^{ae}C^b_{de}+g^{be}C^a_{de})$, from
Killing's equations. Using the relation between $g_{ab}$ and
$h_{ab}$, and the fact that ad is a Lie algebra homomorphism, we find
that
\[
\hat{E}_c(g^{ab})=-A^a_dA^b_e(h^{df}C^e_{cf}+h^{ef}C^d_{cf}).
\]
The expression in the brackets vanishes, as follows easily from
the properties of $h_{ab}$. Thus the generalized Routh equation is
\[
\frac{d}{dt}(g_{ij}v^j)-\onehalf\fpd{g_{jk}}{x^i}v^jv^k=
g_{ij}\left(\dot{v}^j+\conn jklv^kv^l\right)=-\mu_aR^a_{ij}v^j.
\]
But $\mu_a=g_{ab}v^b=\bar{\A}_a^ch_{bc}w^b$; so if we set
$K^a_{ij}=\bar{\A}_b^aR^b_{ij}$, then
$\mu_aR^a_{ij}=h_{bc}K^c_{ij}w^b$. The generalized Routh equation
is therefore equivalent to
\[
\ddot{x}^i+\conn ijk\dot{x}^j\dot{x}^k=g^{im}h_{bc}K^c_{lm}\dot{x}^lw^b.
\]

We also need an equation for $w^a$:\ this comes from the constancy of
$\mu_a$, which we may write as
\[
h_{bc}\frac{d}{dt}(\bar{\A}_a^cw^b)=0.
\]
If we are to understand this equation in the present context, we
evidently need to calculate $\dot{\A}^b_a$.  Now
\[
\dot{\A}^b_a=v^iX_i(\A^b_a)+v^c\tilde{E}_c(\A^b_a)
=v^i\gamma_{ia}^c\A^b_c+v^cC^b_{cd}\A^d_a.
\]
It follows that
\begin{eqnarray*}
h_{bc}\frac{d}{dt}(\bar{\A}_a^c)&=&
-h_{bc}\bar{\A}_a^d\bar{\A}^c_e\dot{\A}^e_d=
-h_{bc}\bar{\A}_a^d\bar{\A}^c_e(v^i\gamma_{id}^f\A^e_f+v^fC^e_{fg}\A^g_d)\\
&=&-h_{bc}\bar{\A}_a^d(v^i\gamma_{id}^c+w^eC^c_{ed}),
\end{eqnarray*}
where in the last step we have again used the fact that ad is a Lie algebra
homomorphism. Now from the skew-symmetry properties of $h_{ab}$ we
obtain
\[
h_{bc}\frac{d}{dt}(\bar{\A}_a^c)=h_{cd}\bar{\A}^d_a(v^i\gamma^c_{ib}+w^eC^c_{eb}),
\]
and therefore
\[
h_{bc}\frac{d}{dt}(\bar{\A}_a^cw^b)=
h_{cd}\bar{\A}^d_a(\dot{w}^c+\gamma^c_{ib}v^iw^b).
\]
The generalized Routh equation and the constancy of momentum together
amount to the mixed first- and second-order equations
\begin{eqnarray*}
\ddot{x}^i+\conn ijk\dot{x}^j\dot{x}^k
&=&g^{im}h_{bc}K^c_{lm}\dot{x}^lw^b\\
\dot{w^a}+\gamma^a_{jb}\dot{x}^jw^b&=&0.
\end{eqnarray*}
These are Wong's equations as they are usually expressed.

\subsection{A Lagrangian with $SE(2)$ as symmetry group}

We now consider the Lagrangian (of simple mechanical type)
\[
L= \onehalf{\dot \z}^2 + \onehalf{\dot \x}^2 +\onehalf {\dot \y}^2
+\onehalf{\dot\theta}^2 + A((\sin\theta){\dot \y} +
(\cos\theta) {\dot \x}){\dot\theta}.
\]
The system is regular if $A^2\neq 1$. The Euler-Lagrange equations are
\[
\ddot \z=0,\quad
\frac{d}{dt} (\dot \x + (A \cos\theta){\dot\theta})=0, \quad
\frac{d}{dt} (\dot \y + (A \sin\theta){\dot\theta})=0, \quad
\ddot\theta + (A\sin\theta) \ddot \y + (A\cos\theta)\ddot \x=0,
\]
and the solution with (for convenience) $\theta_0=0$  is
\begin{eqnarray*}
\lefteqn{(x(t),y(t),z(t),\theta(t))}\\
&&=
\left(\dot{\z}_0t+\z_0,
-A\sin({\dot\theta}_0t)+({\dot \x}_0+A{\dot\theta}_0)t+\x_0,
A\cos({\dot\theta}_0 t) +{\dot \y}_0t +\y_0-A,
{\dot\theta}_0t\right).
\end{eqnarray*}

The system is invariant under the group $SE(2)$, the special Euclidean
group of the plane.  The configuration manifold is ${\bf R}\times
SE(2)$, where $\z$ is the coordinate on ${\bf R}$.  We will use the
trivial connection.  An element of $SE(2)$ can
be represented by the matrix
\[
\left( \begin{array}{ccc} \cos\theta & -\sin\theta & \x \\
\sin\theta &\cos\theta & \y \\ 0&0&1 \end{array}\right).
\]
The identity of the group is $(\x=0,\y=0,\theta=0)$ and the
multiplication is given by
\[
(\x_1,\y_1,\theta_1) * (\x_2,\y_2,\theta_2) =
(\x_2\cos\theta_1-\y_2\sin\theta_1 + \x_1,
\x_2\sin\theta_1+\y_2\cos\theta_1 + \y_1, \theta_1+\theta_2).
\]
The matrices
\[
e_1=\left( \begin{array}{ccc} 0 & 0 & 1 \\
0 &0& 0 \\ 0&0&0\end{array}\right), \qquad
e_2=\left( \begin{array}{ccc} 0 & 0 & 0 \\
0 &0& 1 \\ 0&0&0\end{array}\right), \qquad
e_3=\left( \begin{array}{ccc} 0 & -1 & 0 \\
1 &0& 0 \\ 0&0&0\end{array}\right),
\]
form a basis for the Lie algebra, for which $[e_1,e_2]=0$,
$[e_1,e_3]=e_2$ and $[e_2,e_3]=-e_1$.  The corresponding basis for the
fundamental vector fields is
\[
{\tilde e}_1 = \fpd{}{\x}, \qquad {\tilde e}_2 = \fpd{}{\y}, \qquad
{\tilde e}_3 = - \y\fpd{}{\x}+\x\fpd{}{\y} +\fpd{}{\theta},
\]
and for the invariant vector fields we get
\[
{\hat e}_1 = \cos\theta\fpd{}{\x}+\sin\theta\fpd{}{\y}, \qquad {\hat
e}_2 = -\sin\theta\fpd{}{\x}+\cos\theta \fpd{}{\y} , \qquad {\hat
e}_3 = \fpd{}{\theta}.
\]
One can easily verify that the Lagrangian is invariant.

Before we calculate an expression for the level sets $p_a=\mu_a$, we
will examine the isotropy algebra $\g_\mu$ of a generic point
$\mu=\mu_1e^1+\mu_2e^2+\mu_3e^3$ in $\g^*$. The relations that
characterize an element $\xi=\xi^1e_1+\xi^2e_2+\xi^3e_3$ of $\g_\mu$
are
\[
\xi^3\mu_2 =0, \qquad \xi^3\mu_1 =0, \qquad \xi^1\mu_2-\xi^2\mu_1
=0.
\]
So if we suppose that $\mu_1$ and $\mu_2$ do not both vanish ---
we will take them from now on to be $1$ and $\mu$
respectively --- then a typical element of $\g_\mu$ is $\xi=
\xi^1(e_1+\mu e_2)$. We will also set $\mu_3=0$ for convenience.
Since $\g_\mu$ is 1-dimensional it is of course
Abelian.

Before writing down the coordinate version of the reduced equations in
the previous sections we made two assumptions.  First, we supposed
that a part of the basis of $\g$ was in fact a basis of $\g_\mu$.  So
from now on we will work with a new basis $\{ E_1= e_1+\mu e_2,
E_2=e_2, E_3=e_3\}$, with corresponding notations for the fundamental
and invariant vector fields.  The Lie algebra brackets in this basis
are $[E_1,E_2]=0$, $[E_1,E_3]=-\mu E_1 + (1+\mu^2)E_2$ and
$[E_2,E_3]=-E_1+\mu E_2$.  The momentum vector with which we are
working takes the form
$(1+\mu^2)E^1+\mu E^2$ (with $\mu_3=0$), when written with respect
to the new dual basis.

The second assumption is that we use coordinates
$(\theta^a)=(\theta^A,\theta^\alpha)$ on $G$ such that the fibres
$G\to G/G_\mu$ are given by
$\theta^\alpha=\mbox{constant}$. Then fundamental vector fields for the $G_\mu$-action
on $G$ are of the form $K_A^B\partial/\partial\theta^B$.  The main
advantage of this assumption is that in these coordinates the
expressions in the reduced equations became independent of the
coordinates $\theta^A$. This assumption is not yet satisfied in our
case for the coordinates $(\x,\y,\theta)$. The action of $G_\mu$ on
$G$ is given by the restriction of the multiplication, i.e.\ by
\[
(\x_1)* (\x_2,\y_2,\theta_2) = (\x_2+ \x_1,  \y_2, \theta_2).
\]
We have only one coordinate on $G_\mu$, say $\x'$.
The fundamental vector fields that correspond to this action should
be of the form $K\partial/\partial \x'$. However, in the new basis, vectors in
$\g_\mu$ are of the form $K E_1$, with corresponding fundamental
vector fields
\[
K {\tilde E}_1=K \left(\fpd{}{\x}+\mu\fpd{}{\y}\right).
\]
So we should make a coordinate change $(\x,\y,\theta)\to
(\x',\y',\theta')$, such that
\[
\fpd{}{\x'}=\fpd{}{\x}+\mu\fpd{}{\y}.
\]
This can be done by putting
\[
\x'=\x, \qquad \y'= \y-\mu \x, \qquad \theta'=\theta.
\]
We will then have coordinates $(\x',\y',\theta', \z, \dot
\z)$ on $N_\mu$, and $(\y',\theta', \z, \dot \z)$ on $N_\mu/G_\mu$.
To save typing, we will use $\x$ and $\theta$ for $\x'$ and
$\theta'$, and only make the distinction between $\y$ and $\y'$.

The first goal is to solve the reduced equations on $N_\mu/G_\mu$. They
are of the form
\[
\left\{\begin{array}{lll} {\ddot x}^i  &=&
\Gamma^i(x^j,\theta^\alpha,\dot{x}^j),\\
{\dot \theta}^\alpha &=& \iota^\beta K^\alpha_\beta - {\dot
x}^i{\Lambda}^\alpha_i,\end{array} \right.
\]
For this example there is only one coordinate $\z$ on ${\bf R}$ (we
are using the trivial connection on $SE(2)\times {\bf R}\to {\bf
R}$), but the coordinates $\theta^\alpha$ on $SE(2)/G_\mu$ are
$(\y',\theta)$. The reduced second-order equation in $\z$  above can
be derived from the Lagrangian equation in $\z$  which is simply
\[
\ddot \z=0.
\]
It is therefore not coupled to the first order equation in
$(\y',\theta)$, and its
solution is $\z(t)={\dot \z}_0 t +\z_0$. For the other equations, we
will work first with the variables $(\x,\y,\theta)$, and only make
the change to the new coordinates at the end.

The matrix $(K^\alpha_\beta)$ in the above expressions is determined
by the relation ${\tilde E}_a=K_a^b\partial/\partial\theta^a$. It is the
lower right (2,2)-matrix of
\[
K=\left( \begin{array}{ccc} 1 & 0 & 0 \\
0 & 1 & 0 \\ -\y & \x+\mu \y&1 \end{array}\right).
\]
With the trivial connection, the equations for the other variables
on $N_\mu/G_\mu$ are therefore of the form
\[
{\dot \y}' = \iota^2  + (\x+\mu\y)\iota^3 ,\qquad {\dot\theta}=\iota^3.
\]
We can find the functions $\iota^a$ by solving the expressions
$p_a=\mu_a$ for $v^a$, with $(\mu_1,\mu_2,\mu_3)\in\g^*$ of the form
$((1+\mu^2),\mu,0)$. We get
\begin{eqnarray*}
&&1+\mu^2 = \dot \x+ \mu\dot \y +( A\cos\theta + A\mu\sin\theta)
\dot\theta,\\
&&\mu = \dot \y + (A\sin\theta)\dot\theta,\\
&&0= (A\cos\theta)\dot \x+ (A\sin\theta) \dot \y +\dot\theta-\y(\dot \x
+ (A\cos\theta)\dot\theta)+\x(\dot \y+(A\sin\theta)\dot\theta).
\end{eqnarray*}
At $t=0$, the above equations relate the integration constants and
$\mu$. We will set from now on ${\dot\x}_0=1-A{\dot\theta}_0$,
${\dot \y}_0=\mu$ and $\y_0=\mu \x_0+A{\dot \x}_0+{\dot\theta}_0$.
It is easy to see that the coordinates
$v^a$ with respect to the basis $\{{\tilde E}_1\}$ (with the trivial
connection) are given by
\[
v^1= \dot \x + \y \dot\theta, \quad
v^2 = \dot \y - \mu \dot \x -\mu\y\dot\theta-\x\dot\theta,\quad
v^3=\dot\theta.
\]
After substituting this into the equations for the level set, we
obtain the expressions $v^a=\iota^a$ as functions of
$(\x,\y,\theta)$. After some calculation, the reduced equations
become
\begin{eqnarray*}
{\dot \y}' &=& \frac{A}{A^2-1}\Big(
(\y-\mu\x)(\sin\theta-\mu\cos\theta)-A
(1-\mu^2)\sin\theta\cos\theta-\mu A+2\mu A (\cos\theta)^2\Big),\\
{\dot\theta} &=&  \frac{1}{A^2-1}\Big( \mu\x -\y +
A\cos\theta+A\mu\sin\theta\Big).
\end{eqnarray*}
Observe that we can now replace $(\y-\mu \x)$ everywhere by the new
coordinate $\y'$, so that indeed the $G_\mu$-coordinate $\x'$ does
not appear in the reduced equations. One can verify that the
solution of the above equations, with the integration constants
determined by $\mu$, is
\[
(\y'(t),\theta(t)) = \left(A\cos({\dot\theta}_0t)+A\mu\sin({\dot\theta}_0t)
+(1-A^2){\dot\theta}_0,{\dot\theta}_0t\right).
\]

We will now use the mechanical connection to reconstruct the
$G_\mu$-part $\x(t)$ of the solution. The Hessian of the Lagrangian,
in the basis $\{X=\partial/\partial\z,{\tilde E}_a\}$ is
\[
\left(
\begin{array}{cccc}
1+\mu^2 &\mu&A\cos\theta+A\mu\sin\theta-\y+\mu \x&0\\
\mu&1&A\sin\theta+\x&0\\
A\cos\theta+A\mu\sin\theta-\y+\mu \x& A\sin\theta+\x &
1-2A\y\cos\theta+ 2A\x\sin\theta +\x^2+\y^2 &0\\
0&0&0&1
\end{array}\right).
\]
The determinant of the matrix is $1-A^2$. The vector field
$\clift{\bar X}=\clift X=\partial/\partial\z$ is tangent to the
level sets and horizontal with respect to the mechanical connection
$\mech\Omega$.

In general, we regard the Hessian as a tensor field along the
tangent bundle projection. A basis of vector fields along $\tau$
that lie in the $g$-complement of $\g_\mu$ is
\[
\left\{ {\tilde E}_2 - \frac{\mu}{1+\mu^2}{\tilde E}_1,\,\,
{\tilde E}_3-\frac{1}{1+\mu^2}(A\cos\theta+A\mu\sin\theta-\y'){\tilde E}_1,\,\,
\fpd{}{\z}\right\}.
\]
Notice that they are all basic vector fields along $\tau$ (i.e.\
vector fields on $M$). The reason is that the Lagrangian is of the
simple type. We have seen that in that case the $g$-complement of
$\g_\mu$ defines a connection $\omega^\mu$ on $M \to M/G_\mu$.
The connection tensor $\mech\Omega$ of the mechanical connection is
determined by
\[
\mech\Omega(\clift{\tilde E}_1) = \clift{\tilde E}_1,\qquad
\mech\Omega(\Hb_\alpha) = 0,\qquad
\mech\Omega(\clift{\hat X})=0,\qquad
\mech\Omega(\vlift{\bar X}) =0,
\]
where the vector fields $\{ \hlift{\bar E}_\alpha\}$ that are
horizontal with respect to the mechanical connection $\mech\Omega$ and tangent
to the level set are here
\[
\hlift{\bar E}_2 = \clift{\bar E}_2 -
\frac{\mu}{1+\mu^2}\clift{\tilde E}_1, \qquad
\hlift{\bar E}_3 =\clift{\bar E}_3 -
\frac{1}{1+\mu^2}(A\cos\theta+A\mu\sin\theta-\y')\clift{\tilde E}_1.
\]
The explicit expressions of the $\clift{\bar E}_\alpha$ are not of
direct concern, we only need to know that they are tangent to the
level set and that they differ from $\clift{\hat E}_\alpha$ in a
vertical lift. The vertical part of $\Gamma= \iota^a \clift{\bar
E}_a + v^i \clift{\bar X}_i + \Gamma^i \vlift{\bar X}_i$ (the
restriction of the dynamical vector field to $N_\mu$) is here
\[
\mech\Omega(\Gamma) = \Big(\iota^1 + \frac{\mu}{1+\mu^2}\iota^2 +
\frac{1}{1+\mu^2}(A\cos\theta+A\mu\sin\theta-\y')\iota^3
\Big)\clift{\tilde E}_1.
\]

Before we can write down the explicit form of the reconstruction
equation $g^{-1}\dot{g} = \mech\Omega(\Gamma\circ \hlift{\check v})$,
we need to find the horizontal lift $\hlift{\check v}$ of the reduced
solution $\check v = (\y',\theta,\z,\dot \z)$. It is the curve
$(\x_m,\y',\theta,\z,\dot \z)$ in $N_\mu$ whose tangent vector is
horizontal with respect to the $G_\mu$-mechanical connection. By
construction this means that $\frac{d}{dt} (\tau\circ\hlift{\check v})$
should be $\omega^\mu$-horizontal. If we write in general that
$\frac{d}{dt} (\tau\circ\hlift{\check v}) = v^1 {\tilde E}_1 + v^2
{\tilde E}_2 + v^3 {\tilde E}_3$, then in order for the curve to be
horizontal the $v^a$ must satisfy
\[
v^1 =-v^2\frac{\mu}{1+\mu^2} - v^3
\frac{1}{1+\mu^2}(A\cos\theta+A\mu\sin\theta-\y').
\]
By expressing the $v^a$ as functions of the ${\dot\theta}^a$, we find
that the missing $\hlift\x$ is a solution of
\[
\dot{\hlift\x} = -A{\dot\theta}_0\cos({\dot\theta}_0 t),
\]
from which $\hlift\x(t) = -A\sin({\dot\theta}_0t)+\x_0$. Using this
$\hlift\x$ in the reconstruction equation gives
\[
{\dot \x}_{1} = \iota^1+ \frac{\mu}{1+\mu^2}\iota^2
+\frac{1}{1+\mu^2}(A\cos\theta+A\mu\sin\theta-\y')\iota^3
= 1,
\]
once we have evaluated the functions $\iota^a$ in terms of
$(\hlift\x,\y',\theta,\z,\dot \z)$. So the solution through the
identity is $\x_{1}(t) = t$. The $\x$-part of the complete solution
of the Euler-Lagrange equation is therefore
\[
\x(t) = \x_{1}(t) + \hlift\x(t) = -A\sin({\dot\theta}_0t)+ t+\x_0,
\]
as it be should for the given value of the momentum.

\subsubsection*{Acknowledgements}
The first author is a Guest Professor at Ghent University:\ he is
grateful to the Department of Mathematical Physics and Astronomy at
Ghent for its hospitality.

The second author is currently a Research Fellow at The University
of Michigan through a Marie Curie Fellowship. He is grateful to the
Department of Mathematics for its hospitality. He also acknowledges
a research grant (Krediet aan Navorsers) from the Fund for
Scientific Research - Flanders (FWO-Vlaanderen), where he is an
Honorary Postdoctoral Fellow.

\end{document}